\documentclass[11pt, a4paper, oneside, final, reqno]{amsart}

\usepackage{amssymb}

\newcommand{\calA}{\mathcal{A}}

\newcommand{\calC}{\mathcal{C}}
\newcommand{\calE}{\mathcal{E}}
\newcommand{\calF}{\mathcal{F}}
\newcommand{\calH}{\mathcal{H}}
\newcommand{\calK}{\mathcal{K}}
\newcommand{\calL}{\mathcal{L}}
\newcommand{\calM}{\mathcal{M}}

\newcommand{\calQ}{\mathcal{Q}}

\newcommand{\bbN}{\mathbb{N}}
\newcommand{\bbZ}{\mathbb{Z}}

\newcommand{\bbR}{\mathbb{R}}
\newcommand{\bbC}{\mathbb{C}}
\newcommand{\bbT}{\mathbb{T}}
\newcommand{\bbF}{\mathbb{F}}

\newcommand{\Haus}{\operatorname{Haus}}
\newcommand{\dom}{\operatorname{dom}}

\newcommand{\Hom}{\operatorname{Hom}}

\newcommand{\id}{\operatorname{id}}

\newcommand{\Tr}{\operatorname{Tr}}

\newcommand{\mk}{\operatorname{mk}}

\newcommand{\Span}{\operatorname{span}}

\newcommand{\Cred}{C^*_r}

\newcommand{\ket}[1]{\left|#1\right\rangle}

\theoremstyle{definition}
\newtheorem{definition}{Definition}[section]
\newtheorem{construction}[definition]{Construction}

\theoremstyle{plain}
\newtheorem{theorem}[definition]{Theorem}
\newtheorem{proposition}[definition]{Proposition}
\newtheorem{corollary}[definition]{Corollary}
\newtheorem{lemma}[definition]{Lemma}

\newtheorem{introtheorem}{Theorem}

\theoremstyle{remark}
\newtheorem{remark}[definition]{Remark}
\newtheorem{example}[definition]{Example}

\numberwithin{equation}{section}
\usepackage{hyperref}
\hypersetup{pdftitle={}}
\usepackage[nameinlink, noabbrev, capitalize]{cleveref}
\usepackage[utf8]{inputenc}
\usepackage[backend=biber,style=alphabetic, sorting=nyt, maxnames=50]{biblatex}
\usepackage{xcolor}
\usepackage[all]{xy}

\usepackage[a4paper, total={6.5in, 10in}]{geometry}
\usepackage{microtype}

\usepackage{todonotes}

\title[]{A Quantum Metric Approach to Stabilizer Codes}

\date{}

\author[1]{Are Austad}
\address{Are Austad, Department of Mathematics, University of Oslo, P.O.Box 1053 Blindern, 0316 Oslo, Norway}
\email{areaus@math.uio.no}

\author[2]{Jonas Eidesen}
\address{Jonas Eidesen, School of Mathematics, Statistics and Physics, Herschel Building, Newcastle University, Newcastle upon Tyne, NE1 7RU, UK}
\email{jonas.eidesen@newcastle.ac.uk}



\begin{document}

\begin{abstract}
    We associate a spectral triple with each stabilizer code in a projective error model arising from a nice error basis and equipped with a length function. The spectral triple records the dimension of the code as well as the minimum code distance. For Pauli stabilizer codes of dimension greater than one with Hamming weight, Latrémolière's dual propinquity provides a pseudometric which allows us to compare codes of different physical and logical dimensions. We show that in the resulting topology, minimum distance of a code is locally constant  for qubits and for qudits of dimensions three and five. Thus correction of all errors up to a prescribed weight is a local property in these cases. Moreover, for each prime $p$, we classify the underlying $C^*$-algebras in the metric closure as a subclass of AF algebras, namely those arising as tensor products of $C((\bbZ_p)^r)$, for $r \in \{0,1,2, \ldots , \infty\}$, with matrix algebras of $p$-power order or the UHF algebra of type $p^\infty$. Every such $C^*$-algebra not isomorphic to $M_{p^k}(\bbC)$ for any $k\geq 1$ occurs as the underlying algebra of a point in the metric boundary. In particular, we obtain the CAR algebra as a metric limit of qubit stabilizer codes. 
\end{abstract}

\maketitle

\section{Introduction}
Quantum error correction protects quantum information by encoding it in a subspace $\calC$ of a larger physical Hilbert space $\calH$. The aim is to choose $\calC$ so that a prescribed family of errors can be reversed by a recovery operation. The Knill--Laflamme conditions characterize when this is possible \cite{KnillLaflamme97}: if $P_\calC$ denotes the orthogonal projection onto $\calC$, they require $P_\calC a^* b P_\calC \in \bbC P_\calC$ for all $a,b$ in the prescribed linear space of errors. Stabilizer codes and their systematic study by Gottesman \cite{Gottesman97} provide an algebraic approach to this problem. The code is then a joint eigenspace of commuting error operators, and its logical operations are described by the operators preserving that eigenspace. 

We work in the framework of projective error models (see \cite{Eidesen2026}) associated with nice error bases (see \cite{KnillI96,KnillII96}). Such a model consists of a finite group $G$ and a projectively faithful, irreducible projective representation $\pi$ on $\calH$, with $|G| = \dim(\calH)^2$. The operators $\pi(g)$ form an orthogonal basis of $B(\calH)$ with respect to the Hilbert-Schmidt inner product. In particular, this includes the Pauli error models on qubits and prime-dimensional qudits, which are of particular interest in this article. 
To measure the cost of errors, we equip $G$ with a length function $\ell$, obtained for instance from a cost function on a generating set.
The main example to keep in mind is the Hamming weight, which counts the physical sites on which an operator acts nontrivially. We define the triple $(G,\pi,\ell)$ to be a \emph{projective error model with a weight on $\calH$}, cf. Definition \ref{def:projective error model with a weight}.

These data give a natural connection with noncommutative geometry through spectral triples. Indeed, denoting by $S(\calC)$ and $L(\calC)$ the stabilizer and logical groups of the quantum code $\calC$, respectively, we put $\calL_\calC = L(\calC)/S(\calC)$. After choosing unitary representatives, the projective action on $\calC$ determines a cocycle $\sigma_\calC$, and the logical algebra can be realized as the twisted group $C^*$-algebra $A_\calC = C^*_r(\calL_\calC , \sigma_\calC)$. In fact by Proposition \ref{prop:logical algebra isomorphic to B(C)} we have that $A_\calC \cong B(\calC)$, which recovers the definition of the logical algebra in the operator algebraic quantum error correction framework \cite{KribsLaflammePoulin2005,KribsLaflammePoulin2006,BenyKempfKribs2007I,BenyKempfKribs2007II}. Moreover, any length function $\ell$ on $G$ yields a length function on $\calL_\calC$ through the equation
\begin{align*}
    \ell_\calC (x S(\calC)) = \min_{s \in S(\calC)} \ell(xs).
\end{align*}
$A_\calC$ acts naturally on $\ell^2(\calL_\calC)$ through its left regular representation and $\ell_\calC$ induces a Dirac operator on $\ell^2(\calL_\calC)$ through $D_\calC \delta_x = \ell_\calC(x) \delta_x$. This results in a spectral triple $(A_\calC, \ell^2(\calL_\calC), D_\calC)$ which encodes the algebra of logical observables and the least physical weight of each logical class. In particular, it determines the minimum code distance, as the following theorem makes precise. 

\begin{introtheorem}[cf. Theorem \ref{thm:minimum code distance is minimum of basis logicals}]
    Let $(G, \pi, \ell)$ be a projective error model with a weight on $\calH$ satisfying $|G| = \dim (\calH)^2$. If $\calC$ is a stabilizer code satisfying $\dim(\calC) > 1$, then the minimum code distance is realized as
    \begin{align*}
        \Delta_\calC = \min_{x \in \calL_\calC \setminus \{e\}} \ell_\calC(x) = \inf \{ \Vert [D_\calC, a] \Vert \mid a \in A_\calC, \tau(a) = 0, \tau(a^*a) = 1 \},
    \end{align*}
    where $\tau$ is the normalized trace. In particular, $\Delta_\calC$ is the smallest positive eigenvalue of $D_\calC$. For Pauli codes with Hamming weight, $\Delta_\calC$ is the usual minimum code distance.
\end{introtheorem}

For Pauli codes equipped with Hamming weight and integers $t \geq 0$, the Knill--Laflamme conditions show that correction of the entire space of errors of weight at most $t$ is equivalent to $\Delta_\calC > 2t$. These correction thresholds can therefore be read from the spectral triple.

The passage from spectral triples to metric geometry goes back to Connes \cite{Connes1989, ConnesNCGBook}. A spectral triple determines a commutator seminorm, and the corresponding dual formula measures distances between states. Rieffel developed an abstract framework for Connes's construction by requiring the metric on the state space to induce the weak$^*$-topology \cite{RieffelMetricsCompact98, RieffelMetrics99}. The resulting compact quantum metric spaces admit a generalized notion of Gromov--Hausdorff distance, allowing one to study convergence and continuous families of noncommutative spaces \cite{RieffelqGH04}. 

For our purposes it is important that the metric structure also encodes the underlying $C^*$-algebraic structure. We therefore use Latrémolière's dual propinquity $\Lambda^{\rm dual}$ \cite{Latremoliere2015Dual, Latremoliere2017TriangleIneq}. This yields a complete metric on the class of Leibniz quantum compact metric spaces, up to full quantum isometry. A full quantum isometry is a $*$-isomorphism preserving the Lip-norms. The Leibniz condition is particularly natural for the seminorms considered in this article, as they arise from commutators with Dirac operators. 
Specifically, fix a prime $p$, and let $K_\calC$ be the self-adjoint restriction of the commutator seminorm of a Pauli code with $\dim \calC > 1$ equipped with Hamming weight. We define
\begin{align*}
    d_p (\calC, \calC') = \Lambda^{\rm dual} ((A_\calC, K_\calC), (A_{\calC'}, K_{\calC'})).
\end{align*}
This is a pseudometric on codes, allowing comparisons across different physical and logical dimensions. Its metric quotient will be denoted by $\calQ_p$. Note that by passing to this quotient we forget the physical error labels, but retain the quantum geometry of the logical observables. 
For $p \in \{2,3,5\}$ we in particular prove that the minimum distance is locally constant in this topology. Consequently, for each integer $t \geq 0$, the set of points of $\calQ_p$ correcting all errors of weight at most $t$ is both open and closed, see Theorem \ref{thm:local-error-correction}. Exact correction of bounded weight errors is therefore a local property, without fixing either the physical or logical dimension. 

From the metric on $\calQ_p$ induced by $\Lambda^{\rm dual}$ we also obtain a metric boundary $\partial \calQ_p = \overline{\calQ_p} \setminus \calQ_p$. Here the quantum metric geometry of AF algebras becomes relevant, as they provide natural candidates for limits for our finite-dimensional codes. An early contribution to the study of quantum metric geometry of AF algebras can be found in \cite{ChristensenIvanAF06}. Over the years, AF algebras have been further studied in the setting of quantum compact metric spaces and propinquity, see e.g. \cite{AguilarLatremoliere15, AguilarFelltopAF19, AguilarIndLim21, AAAKMqGH23, AGRKLstrongLeibniz24, AguilarYu24}, though this is not a comprehensive list. 
Our boundary problem starts with a convergent sequence of codes, with no prescribed connecting embeddings between their logical algebras. We show that the Pauli structure nevertheless forces the limiting algebra to be AF, and identify exactly which underlying $C^*$-algebras can occur. Indeed, put 
\begin{align*}
    B_{p,s}=M_{p^s}(\bbC)\quad(s<\infty),\qquad  B_{p,\infty}=\bigotimes_{j\geq1}M_p(\bbC).
\end{align*}

\begin{introtheorem}[cf. Theorems \ref{thm:limits-are-AF} and \ref{thm:AF-boundary-realization}]
    The $C^*$-algebras underlying points of $\overline{\calQ_p}$ are exactly 
    \begin{align*}
        C((\bbZ_p)^r)\otimes B_{p,s}, \qquad r,s\in\{0,1,2,\ldots,\infty\},
    \end{align*}
    where $(\bbZ_p)^{\infty}$ carries the product topology. Every listed algebra admits a metric spectral triple whose associated quantum metric space belongs to $\overline{\calQ_p}$. Every algebra underlying a boundary point is therefore AF, and each listed algebra not isomorphic to $M_{p^k} (\bbC)$ for any $k \geq 1$ occurs at a boundary point. 
\end{introtheorem}
The realizations use ordinary physical Hamming weights throughout. 
In particular, the UHF algebra of type $p^\infty$ occurs as the underlying algebra of a boundary point. For $p=2$, a concrete sequence converging to the CAR algebra is obtained by taking tensor products of successively concatenated five-qubit codes. The $m$th code is
$\bigotimes_{j=1}^m\calF^{(j)}$, where $\calF^{(j)}$ is the $j$-fold concatenation of the $[[5,1,3]]_2$-code from Lemma \ref{lemma:513-protector-code}, see  Theorem \ref{theorem:CAR-limit}.

Section \ref{sec:preliminaries} recalls the required background on projective error models, stabilizer codes, spectral triples and quantum compact metric spaces. Section \ref{sec:quantum-codes-as-spectral-triples} constructs the logical spectral triple and recovers the minimum code distance. Section \ref{sec:quantum-distances-between-codes} studies the resulting topology, its metric limits, and local error correction. 

\subsubsection*{Acknowledgments} The first author wishes to thank NTNU Trondheim for its hospitality during a research stay during fall 2026, during which parts of this research were conducted. 

The first author was funded by the Research Council of Norway projects 324944 and 359602.
The second author was supported by the Engineering and Physical Sciences Research Council [grant number UKRI3684].

\section{Preliminaries}
\label{sec:preliminaries}
\subsection{Notation and conventions}

Any Hilbert space considered will be complex and nonzero. Unless otherwise specified, groups are discrete. We denote the identity of a group $G$ by $e$, or by $e_G$, when the group needs to be specified.
Similarly, if $\calH$ is a Hilbert space, we will denote the identity operator on $\calH$ by $I_\calH$. In the specific case where $\calH = \bbC^n$ we will use the notation $I_n := I_\calH$. 
The bounded linear operators on a Hilbert space $\calH$ will be denoted by $B(\calH)$, and in the case where $\calH = \bbC^n$ we make the identification $B(\calH) = M_n(\bbC)$, where $M_n(\bbC)$ denotes the space of all $n \times n$ complex matrices. 
For a closed subspace $\calC \subset \calH$, we will use the notation $P_\calC \in B(\calH)$ to denote the orthogonal projection onto the subspace $\calC$. 
For finite-dimensional $\calH$, we write $\Tr$ for the unnormalized trace on $B(\calH)$. If we need to emphasize that it is the trace on $B(\calH)$, we will use the notation $\Tr_\calH$, if $\calH = \bbC^n$ we will instead use the notation $\Tr_n$.

We will use the following notation for cyclic groups of order $n$:
\begin{equation*}
    \bbZ_n = \{0,1,2,\dots,n-1\},
\end{equation*}
where the group operation on $\bbZ_n$ is addition modulo $n$. The notation $\bbT$ will be used for the circle group, explicitly
\begin{equation*}
    \bbT = \{ z \in \bbC \mid |z| = 1 \},
\end{equation*}
and the group operation on $\bbT$ is multiplication. The group of unitary operators on $\calH$ will be denoted by $U(\calH)$. The projective unitary group on $\calH$ is then defined to be the following quotient group:
\begin{equation*}
    PU(\calH) := U(\calH)/\{ z I_\calH \mid z \in \bbT \}.
\end{equation*}
We will always let $q \colon U(\calH) \to PU(\calH)$ denote the quotient map.

We will furthermore make use of the following notation: We will let $\{\ket{i}\}_{i = 0}^{d - 1}$ denote the standard basis for $\bbC^d$, and use the shorthand
\begin{equation*}
    \ket{j_1j_2 ...j_n} := \ket{j_1} \otimes \ket{j_2} \otimes \cdots \otimes \ket{j_n}
\end{equation*}
to denote the basis elements of $\bbC^{d_1} \otimes \bbC^{d_2} \otimes \cdots \otimes \bbC^{d_n}$, where $0 \leq j_1 \leq (d_1 - 1)$, $0 \leq j_2 \leq (d_2 - 1)$, and so on.

\subsection{Projective error models}\label{sec:projective error models}

We will employ the definition of \emph{projective error models} as found in \cite{Eidesen2026} to highlight the group-theoretic and representation-theoretic aspects of our theory. Before we state the definition of a projective error model we need to establish some preliminaries on projective representation theory. For a textbook recounting of the basics of projective representation theory see for example \cite[Chapters 7-9]{CeccheriniSilbersteinTullio22}.

\begin{definition}\label{def:projective representation}
    Let $G$ be a group and $\calH$ be a Hilbert space. A function $\pi \colon G \to U(\calH)$ is a \emph{projective representation of $G$ on $\calH$} if $\pi(e) = I_\calH$ and there exists a function $\sigma_\pi \colon G \times G \to \bbT$ such that
    \begin{equation*}
        \pi(x)\pi(y) = \sigma_\pi(x,y)\pi(xy)
    \end{equation*}
    for all $x,y \in G$. Note that if such a function $\sigma_\pi$ exists, then it is uniquely determined by the function $\pi$.
\end{definition}

The existence of the function $\sigma_\pi$ is equivalent to asking that the composition $q \circ \pi \colon G \to PU(\calH)$ is a group homomorphism. This perspective is useful for the following definition.

\begin{definition}\label{def:projectively faithful}
    Let $G$ be a group and $\pi$ be a projective representation of $G$ on $\calH$. We say that $\pi$ is \emph{projectively faithful} if the group homomorphism $q \circ \pi \colon G \to PU(\calH)$ is injective.
\end{definition}

By associativity of the multiplication in $U(\calH)$, the function $\sigma_\pi$ will satisfy the equations
\begin{equation}\label{eqn:normalized 2-cocycle equation}
    \sigma_\pi(x,y)\sigma_\pi(xy,z) = \sigma_\pi(x,yz)\sigma_\pi(y,z) \quad \text{and} \quad \sigma_\pi(x,e) = 1 = \sigma_\pi(e,x)
\end{equation}
for all $x,y,z \in G$. The first of these equations is called the \emph{$2$-cocycle equation}, and it motivates the following standard definition.

\begin{definition}\label{def:2-cocycle}
    Let $G$ be a group. A \emph{(normalized) $2$-cocycle on $G$ (with values in $\bbT$)} is a function $\sigma \colon G \times G \to \bbT$ satisfying \eqref{eqn:normalized 2-cocycle equation}. The collection of all normalized $2$-cocycles on $G$ is denoted by $Z^2(G;\bbT)$. If there exists a function $\zeta \colon G \to \bbT$ such that
    \begin{equation*}
        \sigma(x,y) = \zeta(x) \, \zeta(y) \, \overline{\zeta(xy)}
    \end{equation*}
    for all $x,y \in G$, then we say that $\sigma$ is a \emph{$2$-coboundary on $G$}. The collection of all $2$-coboundaries on $G$ is denoted by $B^2(G;\bbT)$. Point-wise multiplication endows $Z^2(G;\bbT)$ with the structure of an abelian group, and we identify $B^2(G;\bbT)$ as a subgroup of $Z^2(G;\bbT)$. The quotient group
    \begin{equation*}
        H^2(G;\bbT) := Z^2(G;\bbT)/B^2(G;\bbT)
    \end{equation*}
    is called the \emph{$2$nd cohomology group of $G$ (with values in $\bbT$)}.
\end{definition}

We will not need any of the theory of group cohomology for this paper, but the interested reader may find \cite{Brown82} a valuable resource for learning more about the cohomology of groups.

A standard approach to studying projective representations of a group $G$ is to first fix a $2$-cocycle $\sigma$ on $G$ and look at all projective representations $\pi$ of $G$ that satisfy $\sigma_\pi = \sigma$. In this setting, we will use the term \emph{$\sigma$-projective representation} to emphasize the dependence on the $2$-cocycle $\sigma$. The most important example of a $\sigma$-projective representation is the following.

\begin{example}\label{ex:left regular representation}
    Let $G$ be a group and $\sigma \in Z^2(G;\bbT)$. The \emph{$\sigma$-twisted left regular representation of $G$} is the function $\lambda^\sigma \colon G \to U(\ell^2(G))$ defined by linearly extending the following equation for $x,g \in G$:
    \begin{equation*}
        \lambda^\sigma(x)\delta_g = \sigma(x,g)\delta_{xg}.
    \end{equation*}
    Here $\delta_g$ is the characteristic function on the singleton $\{g\}$. We will also use the notation $\lambda^\sigma_x := \lambda^{\sigma}(x)$.
\end{example}

\begin{definition}\label{def:intertwiner}
    Let $G$ be a group, $\sigma \in Z^2(G;\bbT)$, and $\pi_1$, $\pi_2$ be two $\sigma$-projective representations of $G$ on $\calH_1$, $\calH_2$ respectively. An \emph{intertwiner} from $\pi_1$ to $\pi_2$ is a bounded linear map $T \colon \calH_1 \to \calH_2$ such that
    \begin{equation*}
        T \pi_1(x) = \pi_2(x) T
    \end{equation*}
    for all $x \in G$. The collection of all intertwiners from $\pi_1$ to $\pi_2$ is denoted by $\Hom_G(\pi_1,\pi_2)$. If $T$ is furthermore a bounded linear isomorphism between $\calH_1$ and $\calH_2$, then we say that $\pi_1$ and $\pi_2$ are \emph{isomorphic}, and we will use the notation $\pi_1 \cong \pi_2$ to denote that such an intertwiner exists.
\end{definition}

Note that this definition endows the collection of $\sigma$-projective representations of $G$ with the structure of a category. We will not make use of this fact.

\begin{definition}\label{def:invariant subspace and irreducible projective representation}
    Let $G$ be a group and $\pi$ be a projective representation of $G$ on $\calH$. A closed subspace $\calK \subset \calH$ is \emph{$G$-invariant} if $\pi(x)\ket{\psi} \in \calK$ for all $x \in G$ and all $\ket{\psi} \in \calK$. The projective representation $\pi$ is said to be \emph{irreducible} if the only $G$-invariant subspaces of $\calH$ are $\calH$ and $\{0\}$.
\end{definition}

In the case where $G$ is a finite group, we have that analogues of Schur's lemma and Maschke's theorem hold for $\sigma$-projective representations of $G$. As we will not directly have a use for these results, we will not state them here. We will however have a use for the definition of the character of a projective representation, cf. \cite{Cheng15} for an introduction on the character theory of projective representations. We state the main definition here:

\begin{definition}\label{def:character}
    Let $G$ be a group and $\pi$ be a projective representation of $G$ on a finite-dimensional Hilbert space $\calH$. The \emph{character of $\pi$} is the function $\chi_\pi \colon G \to \bbC$ defined by setting
    \begin{equation*}
        \chi_\pi(x) = \Tr(\pi(x)).
    \end{equation*}
\end{definition}

With this we have everything we need to give a definition of a projective error model.

\begin{definition}\label{def:projective error model}
    Let $\calH$ be a finite dimensional Hilbert space. A \emph{projective error model on $\calH$} is a pair $(G,\pi)$ where $G$ is a finite group, and $\pi$ is a projectively faithful, irreducible projective representation of $G$ on $\calH$.
\end{definition}

The most important example of a projective error model is the following, giving rise to the Pauli projective error model, which will serve as a guiding example throughout this paper.

\begin{example}\label{ex:projective Pauli error model}
    The Pauli projective error model is given by $(P_n,\pi_n)$ where $P_n = (\bbZ_2 \times \bbZ_2)^{n}$ and $\pi_n$ is the $n$-fold product representation of the only projectively faithful, irreducible projective representation of $\bbZ_2 \times \bbZ_2$. Explicitly, if
    \begin{equation*}
        X =
        \begin{bmatrix}
            0 & 1 \\
            1 & 0
        \end{bmatrix}, \ 
        Z =
        \begin{bmatrix}
            1 & 0 \\
            0 & -1 
        \end{bmatrix},
    \end{equation*}
    we define $\pi_n \colon P_n \to U((\bbC^2)^{\otimes n})$ by
    \begin{equation*}
        \pi_n (a_1,b_1,\dots,a_n,b_n) = X^{a_1} Z^{b_1} \otimes \cdots \otimes X^{a_n}Z^{b_n}.
    \end{equation*}
\end{example}

We will mainly restrict ourselves to the study of projective error models $(G,\pi)$ on $\calH$ where $|G| = \dim(\calH)^2$. This setting was first studied by Knill, cf. \cite{KnillI96,KnillII96} where the term \emph{nice error basis} is used for the image $\pi(G) \subset B(\calH)$. Note that in this context, a group $G$ that admits an irreducible projective representation $\pi \colon G \to U(\calH)$ with the property that $|G| = \dim(\calH)^2$ is also referred to as an \emph{index group} in \cite{KlappeneckerRottelerI02,KlappeneckerRottelerII02}, while in the mathematics literature such a group $G$ is referred to as a \emph{central type factor group}.

If one does not include the assumption that $|G| = \dim(\calH)^2$, then a projective error model is the same notion as that of a \emph{nice error frame} introduced in \cite{ChienWaldron17}. For a thorough account of nice error frames, one can consult the book \cite{Waldron18}.

\subsection{Stabilizer codes}\label{sec:stabilizer codes}

The most general definition of a quantum code is that a quantum code is nothing more than a subspace $\calC$ of a Hilbert space $\calH$. However, we are mostly interested in constructing such subspaces given some extra (usually algebraic) data. Pauli stabilizer codes, as introduced by Gottesman in \cite{Gottesman97} are an important example. In this section we recall the definition of a stabilizer code as it is defined in terms of a projective error model in \cite{Eidesen2026} as well as some of the main results we will need about these stabilizer codes.

\begin{definition}\label{def:stabilier code}
    Let $(G,\pi)$ be a projective error model on $\calH$. A nonzero subspace $\calC \subset \calH$ is a \emph{stabilizer code} if there exists a normal subgroup $S \subset G$ and a function $\lambda \colon S \to \bbT$ such that
    \begin{equation}\label{eqn:stabilizer code}
        \calC = \{ \ket{\psi} \in \calH \mid \pi(s)\ket{\psi} = \lambda(s)\ket{\psi} \text{ for all } s \in S \}.
    \end{equation}
\end{definition}

The following result gives a description of when a normal subgroup $S$ of $G$ admits a function $\lambda \colon S \to \bbT$ such that the space on the right hand side of \eqref{eqn:stabilizer code} is nonzero. See \cite[Theorem 5.5]{Eidesen2026} for a proof.

\begin{proposition}\label{prop:nontriviality of stabilizer code}
    Let $(G,\pi)$ be a projective error model on $\calH$. If $S$ is a normal subgroup of $G$ then there exists a function $\lambda \colon S \to \bbT$ such that
    \begin{equation*}
        \{ \ket{\psi} \in \calH \mid \pi(s)\ket{\psi} = \lambda(s)\ket{\psi} \text{ for all } s \in S \}
    \end{equation*}
    is nonzero if and only if $S$ is abelian and $\sigma_\pi|_{S \times S}$ is a $2$-coboundary.
\end{proposition}

It will be beneficial to define the following two subgroups of $G$ in our study of stabilizer codes.

\begin{definition}\label{def:logical group and stabilizer group}
    Let $(G,\pi)$ be a projective error model on $\calH$ and $\calC \subset \calH$ be a subspace. Define the \emph{logical group of $\calC$ (with respect to $(G,\pi)$)} and the \emph{stabilizer group of $\calC$ (with respect to $(G,\pi)$)} to be the following two subgroups of $G$ respectively (recall that $P_\calC \in B(\calH)$ denotes the orthogonal projection onto $\calC$):
    \begin{align*}
        L_{(G,\pi)}(\calC) & = \{ x \in G \mid \pi(x) P_\calC = P_\calC \pi(x) \}, \\
        S_{(G,\pi)}(\calC) & = \{ x \in G \mid P_\calC \pi(x) P_\calC \in \bbT P_\calC \}.
    \end{align*}
    When the projective error model is fixed, we will use the notation $L(\calC) := L_{(G,\pi)}(\calC)$ and $S(\calC) := S_{(G,\pi)}(\calC)$. Note that $S(\calC)$ is always a subgroup of $L(\calC)$, and by definition of the group $S(\calC)$, there exists a function $\lambda_\calC \colon S(\calC) \to \bbT$ defined by the equation
    \begin{equation*}
        P_\calC \pi(s) P_\calC = \lambda_\calC(s) P_\calC.
    \end{equation*}
\end{definition}

We collect several results from \cite{Eidesen2026} that will be useful for us in the following.

\begin{proposition}\label{prop:collection of stabilizer code results}
    Let $(G,\pi)$ be a projective error model on $\calH$ such that $|G| = \dim(\calH)^2$, and $\calC \subset \calH$ be a stabilizer code. If $S$ is a normal subgroup of $G$ and $\lambda \colon S \to \bbT$ is a function such that \eqref{eqn:stabilizer code} holds, then:
    \begin{enumerate}
        \item $L(\calC) = \{ x \in G \mid \lambda(s) = \sigma_\pi(x^{-1},s) \overline{\sigma_\pi(x^{-1}sx,x^{-1})} \lambda(x^{-1}sx) \text{ for all } s \in S\}$,
        \item $S(\calC) = S$ and $\lambda_\calC = \lambda$, 
        \item $|G| = |L(\calC)| \, |S(\calC)|$,
        \item $\dim(\calC) = \dim(\calH)/|S(\calC)|$,
        \item and finally we have that
        \begin{equation*}
            P_\calC = \frac{1}{|S(\calC)|}\sum_{s \in S(\calC)} \overline{\lambda_\calC(s)}\pi(s).
        \end{equation*}
    \end{enumerate}
\end{proposition}

The quotient group $L(\calC)/S(\calC)$ will feature heavily throughout this paper, we therefore make the following definition for notational convenience.

\begin{definition}\label{def:pure logical group}
    Let $(G,\pi)$ be a projective error model on $\calH$ such that $|G| = \dim(\calH)^2$, and $\calC \subset \calH$ be a stabilizer code. Define $\calL_\calC := L(\calC)/S(\calC)$.
\end{definition}

We will have use for the following result.

\begin{lemma}\label{lemma:size of pure logical group}
    Let $(G,\pi)$ be a projective error model on $\calH$ such that $|G| = \dim(\calH)^2$, and $\calC \subset \calH$ be a stabilizer code. Then $\dim(B(\calC)) = |\calL_\calC|$.
\end{lemma}
\begin{proof}
    By Proposition \ref{prop:collection of stabilizer code results} we compute that
    \begin{equation*}
        \dim(B(\calC))
        = \dim(\calC)^2
        = \frac{\dim(\calH)^2}{|S(\calC)|^2}
        = \frac{|G|}{|S(\calC)|^2}
        = \frac{|L(\calC)||S(\calC)|}{|S(\calC)|^2}
        = \frac{|L(\calC)|}{|S(\calC)|}
        = |\calL_\calC|. \qedhere
    \end{equation*}
\end{proof}

\subsection{Twisted group \texorpdfstring{$C^*$}{C*}-algebras}\label{sec:twisted group C*-algebras}

Fix a group $G$ and a $2$-cocycle $\sigma \in Z^2(G;\bbT)$. Let $C_c(G,\sigma)$ denote the $*$-algebra of complex functions on $G$ with finite support equipped with the following $\sigma$-twisted convolution product and involution:
\begin{equation*}
    (f_1 *_\sigma f_2)(x) = \sum_{y \in G}\sigma(y,y^{-1}x)f_1(y)f_2(y^{-1}x), \quad \text{and} \quad f^*(x) = \overline{\sigma(x,x^{-1})} \, \overline{f(x^{-1})}
\end{equation*}
for $f,f_1,f_2 \in C_c(G,\sigma)$ and $x \in G$. The $\sigma$-twisted left regular representation of $G$ induces a $*$-homomorphism $\lambda^\sigma \colon C_c(G,\sigma) \to B(\ell^2(G))$ defined by
\begin{equation*}
    \lambda^{\sigma}(f)\xi(x) = \sum_{y \in G} \sigma(y,y^{-1}x)f(y)\xi(y^{-1}x),
\end{equation*}
for $f \in C_c(G,\sigma)$, $\xi \in \ell^2(G)$, and $x \in G$. With this definition we see that $\lambda^{\sigma}(\delta_x) = \lambda^\sigma_x$ as defined in Example \ref{ex:left regular representation}. This allows us to define the following.

\begin{definition}\label{def:twisted group C*-algebra}
    The \emph{reduced $\sigma$-twisted group $C^*$-algebra of $G$}, which we denote by $\Cred(G,\sigma)$, is the operator norm closure of the image of $C_c(G,\sigma)$ in $B(\ell^2(G))$ under $\lambda^\sigma$.
\end{definition}

There is furthermore a faithful canonical tracial state $\tau_\sigma \in S(\Cred(G,\sigma))$ determined by
\begin{equation*}
    \tau_\sigma(\lambda^\sigma(f)) = f(e),
\end{equation*}
for $f \in C_c(G,\sigma)$.

If $G$ is amenable (in particular, if $G$ is finite) and $ \pi \colon G \to U(\calH)$ is a $\sigma$-projective representation, then we get a $*$-representation of the reduced $\sigma$-twisted group $C^*$-algebra of $G$, $\pi \colon \Cred(G,\sigma) \to B(\calH)$, which is determined by
\begin{equation*}
    \pi(f) = \sum_{x \in G}f(x)\pi(x),
\end{equation*}
for $f \in C_c(G,\sigma)$. The induced $*$-representation is called the \emph{integrated representation of $\pi$}. It is then straightforward to check that this gives the following result.

\begin{lemma}\label{lemma:correspondence between projective representations and *-representations}
    Let $G$ be an amenable group and $\sigma \in Z^2(G;\bbT)$. Then there is a bijective correspondence between $\sigma$-projective representations of $G$ and unital $*$-representations of $\Cred(G,\sigma)$.
\end{lemma}

\subsection{Length functions and spectral triples}\label{sec:length functions and spectral triples}

Quantum error correction is more than choosing a Hilbert space $\calH$ to model a physical system and then deciding on a suitable subspace $\calC \subset \calH$ to act as a quantum code. Any physical system should also be equipped with a notion of severity for any potential error that can occur. For our purposes, this extra data will be geometric in nature, and often be induced from a length function:

\begin{definition}\label{def:length function}
    Let $G$ be a group. A \emph{length function on $G$} is a function $\ell \colon G \to [0,\infty)$ that satisfies the following:
    \begin{enumerate}
        \item $\ell(x) = 0$ if and only if $x = e$,
        \item $\ell(x^{-1}) = \ell(x)$ for all $x \in G$, and
        \item $\ell(xy) \leq \ell(x) + \ell(y)$ for all $x,y \in G$.
    \end{enumerate}
    We say that the length function $\ell$ is \emph{proper} if it is proper as a continuous map, i.e. the inverse image of a compact subset of $[0,\infty)$ under $\ell$ is finite.
\end{definition}

Groups with proper length functions are a valuable source for noncommutative unital spectral triples. The following definition of a unital spectral triple can be found in for instance \cite{Connes1989, ConnesNCGBook, Kaad24}.

\begin{definition}\label{def:spectral triple}
    A \emph{unital spectral triple} $(A,\calH,D)$ consists of a unital $C^*$-algebra $A$, a faithful unital $*$-representation $\pi \colon A \to B(\calH)$, and a (usually unbounded) self-adjoint operator $D \colon \dom(D) \to \calH$, with $\dom(D) \subset \calH$. These need to satisfy the following:
    \begin{enumerate}
        \item There exists a dense $*$-algebra $\calA \subset A$ containing the unit $1_A$ such that for every $a \in \calA$ the commutator
        \begin{equation*}
            [D,\pi(a)] = D \pi(a) - \pi(a) D
        \end{equation*}
        is defined on $\dom D$ and extends to a bounded operator on $\calH$. In particular, $\pi(a) \dom D \subseteq \dom D$. 
        \item The resolvent $(i + D)^{-1}$ is a compact operator.
    \end{enumerate}
\end{definition}

We will often think of $A$ as a concrete $C^*$-subalgebra of $B(\calH)$, hence we will usually not use any notation for the map $A \hookrightarrow B(\calH)$.

If $G$ is a group, $\sigma \in Z^2(G;\bbT)$, and $\ell \colon G \to [0,\infty)$ is a proper length function, then we attain a spectral triple $(\Cred(G,\sigma),\ell^2(G),D_\ell)$, where $D_\ell$ is the self-adjoint multiplication operator with domain $\{\xi \in \ell^2 (G) \mid \sum_x \ell(x)^2 |\xi(x)|^2 < \infty \}$, defined by
\begin{equation*}
    (D_\ell \xi)(x) = \ell(x)\xi(x),
\end{equation*}
for $\xi \in \dom (D_\ell)$ and $x \in G$. The subspace $C_c(G) \subseteq \ell^2(G)$ is a core for $D_\ell$. Properness of $\ell$ gives compact resolvent, and the triangle inequality gives $\Vert [D_\ell, \lambda_x^\sigma ] \Vert \leq \ell(x)$. Thus $C_c(G, \sigma)$ supplies the required dense $*$-algebra.

To any unital spectral triple one can associate to it a Lipschitz seminorm, which we recall the definition of here.

\begin{definition}\label{def:Lipschitz seminorm}
    Let $A$ be a unital $C^*$-algebra. A seminorm $K \colon A \to [0,\infty]$ is called a \emph{Lipschitz seminorm} if:
    \begin{enumerate}
        \item the set $\dom(K) := \{ a \in A \mid K(a) < \infty \}$ is dense in $A$,
        \item $K(a) = K(a^*)$ for all $a \in A$, and
        \item $\bbC 1_A \subset \ker(K) := \{ a \in A \mid K(a) = 0 \}$.
    \end{enumerate}
    To such a seminorm one can associate to it the \emph{Monge-Kantorovi\v{c} (extended) metric},
    \begin{equation*}
        \mk_K \colon S(A) \times S(A) \to [0,\infty],
    \end{equation*}
    defined by
    \begin{equation*}
        \mk_K(\varphi,\psi) = \sup\{ |\varphi(a) - \psi(a)| \mid a \in A, \ K(a) \leq 1 \}.
    \end{equation*}
    If we furthermore have that $\mk_K$ metrizes the weak$^*$-topology on $S(A)$, then we say that $K$ is a \emph{Lip-norm} and that the pair $(A,K)$ is a \emph{compact quantum metric space}.
\end{definition}

The Lipschitz seminorm associated to the unital spectral triple $(A,\calH,D)$ is defined by setting
\begin{equation*}
    K_D(a) =
    \begin{cases}
        \|[D,a]\| & \text{if } a \dom(D) \subseteq \dom(D) \text{ and } [D,a] \text{ is bounded,} \\
        \infty & \text{otherwise.}
    \end{cases}
\end{equation*}
In the special case where the unital spectral triple is coming from a group $G$ with a proper length function $\ell$, we will denote the associated Lipschitz seminorm by $K_\ell$. Note that in the case where $G$ is finite, $K_\ell$ will always be a Lip-norm.

\subsection{Propinquity}
\label{subsec:propinquity}
The seminorm associated to a spectral triple is defined on arbitrary elements of the $C^*$-algebra, while the propinquity compares seminorms defined on self-adjoint elements. These give the same induced Monge-Kantorovi\v{c} metric on the state spaces \cite[Remark 2.2]{Rieffel2025}. That is, for  $\tilde{K}$ a Lipschitz seminorm as above, and $K = \tilde{K}\vert_{A_{\rm sa}}$ we have
\begin{align*}
    \mk_{\tilde{K}}(\phi, \psi) = \mk_{K}(\phi, \psi) := \sup\{ \vert \phi(a) - \psi(a) \vert \mid a\in A_{\rm sa}, K(a) \leq 1 \}
\end{align*}
for $\phi, \psi \in S(A)$. As is customary in the field, we denote the finite domain of $K$ by $\dom(K)$ and extend $K$ by infinity  on $A_{\rm sa} \setminus \dom (K)$. 

A \emph{quantum compact metric space} is a unital $C^*$-algebra $A$ with a seminorm $K$ on a norm-dense subspace $\dom(K)\subseteq A_{\mathrm{sa}}$ containing $1_A$. Its kernel must be $\bbR 1_A$, and $\mk_K$ must induce the weak$^*$-topology on $S(A)$. 
Define
\begin{align*}
    a \circ b = \frac{1}{2} (ab + ba), \qquad \{ a,b\} = \frac{1}{2i} (ab - ba), \qquad \text{for } a,b\in A.
\end{align*}
$K$ is a \emph{Leibniz seminorm} if for all $a,b \in \dom K$
\begin{align}\label{eq:Leibniz-condition}
    K(a\circ b) \leq \Vert a \Vert K(b) + \Vert b \Vert K(a), \qquad K(\{a,b\}) \leq \Vert a \Vert K(b) + \Vert b \Vert K(a). 
\end{align}
A \emph{Leibniz quantum compact metric space} is a pair $(A,K)$ consisting of a quantum compact metric space whose seminorm $K$ is a Leibniz seminorm which is lower semicontinuous with respect to the norm of $A$. 

For a densely defined seminorm with scalar kernel, the following compactness criterion gives a useful test for whether the Monge-Kantorovi\v{c} metric induces the weak$^*$-topology: for one, or equivalently any, $\mu \in S(A)$, the set
\begin{align}\label{eq:Rieffel-compactness-criterion}
    E_{K, \mu} := \{a \in A_{\rm sa} \mid K(a) \leq 1, \mu(a) = 0 \} 
\end{align}
must be norm totally bounded \cite[Theorem 2.10]{Latremoliere2015Dual}.

Note that for a seminorm induced by a spectral triple, the identity $[D,ab] = [D,a]b + a[D,b]$ immediately implies \eqref{eq:Leibniz-condition}. We take the self-adjoint realization of $D$, and for $a = a^*$, $K(a)=\|[D,a]\|$ if $a$ preserves $\dom(D)$ and $[D,a]$ is bounded there, and $K(a)=\infty$ otherwise. This convention gives lower semicontinuity of $K$. 

For a finite twisted group algebra with length $\ell$, we may write $a = \sum_g a_g \lambda^\sigma_g$, giving $\Vert [D_\ell, a ]\delta_e \Vert^2 =\sum_g \vert a_g \vert ^2 \ell(g)^2$, which immediately implies $\ker K = \bbR \cdot 1$. Finite-dimensionality then gives \eqref{eq:Rieffel-compactness-criterion}. 

We now turn to defining distances between quantum compact metric spaces. Let us be given Leibniz quantum compact metric spaces $(A_1, K_1)$ and $(A_2,K_2)$. A \emph{tunnel} from $(A_1,K_1) $ to $(A_2, K_2)$ consists of a Leibniz quantum compact metric space $(C,M)$  and surjective unital $*$-homomorphisms $\rho_1\colon C \to A_1$ and $\rho_2 \colon C \to A_2$ for which the quotient seminorms agree with $K_1$ and $K_2$, that is, for $i =1,2$
\begin{align*}
    K_i (a) = \inf \{ M(x) \mid x \in C_{\rm sa}, \rho_i (x) = a \}, \qquad a\in (A_i)_{\rm sa}. 
\end{align*}
As such, the map $\rho_i^* \colon \phi \mapsto \phi \circ \rho_i $ is an isometric embedding of $(S(A_i), \mk_{K_i})$ into $(S(C), \mk_M)$ for $i=1,2$. We will refer to $(A_i, K_i)$, $i=1,2$ as the endpoints of the tunnel. 
Denoting by $\Haus_d (X,Y)$ the Hausdorff distance between two subsets $X,Y$ of a metric space $(Z,d)$, we then define the \emph{extent} of the tunnel $\mathfrak{t} =(C,M , \rho_1, \rho_2)$ by
\begin{align}\label{eq:def-extent}
    \chi (\mathfrak{t}) = \max_{i\in \{1,2\}} \Haus_{\mk_M} (S(C), \rho_i^*(S(A_i))). 
\end{align}
The \emph{dual propinquity} of $(A_1, K_1)$ and $(A_2, K_2)$ is then
\begin{align}\label{eq:def-dual-propinquity}
    \Lambda^{\rm dual}((A_1,K_1), (A_2, K_2)) := \inf \{ \chi(\mathfrak{t}) \mid \mathfrak{t}\text{ is a tunnel from $(A_1, K_1)$ to $(A_2, K_2)$} \}. 
\end{align}
This distance satisfies the triangle inequality, and is complete up to full quantum isometry, that is, two pairs have distance zero exactly when there is a $*$-isomorphism $\theta \colon A_1 \to A_2$ such that $K_2 (\theta(a)) = K_1(a)$ on $(A_1)_{\rm sa}$, see \cite[Theorem 3.8]{Latremoliere2017TriangleIneq}. 

We finish this section by noting two consequences of \eqref{eq:def-extent}. First, if $\chi(\mathfrak{t}) \leq \varepsilon$, then every state on $C$ is within $\varepsilon$ of a pulled back state from either endpoint. Expressing the norm of a self-adjoint element $x \in \dom (M)$ through the supremum over absolute values of states on $C$, we find
\begin{align}\label{eq:lift-norm-bounds}
    \Vert \rho_i (x) \Vert \leq \Vert x \Vert \leq \Vert \rho_i (x) \Vert + \varepsilon M(x)
\end{align}
for all $x \in \dom (M)$, and for $i=1,2$. 

For the convergence results in Section \ref{sec:quantum-distances-between-codes}, conditional expectations give a direct way to bound the dual propinquity between two Leibniz quantum compact metric spaces. The following result will be useful. 

\begin{lemma}\label{lemma:cp-approximation-map}
    Let $(A, K_A)$ and $(B, K_B)$ be Leibniz quantum compact metric spaces. Suppose $B \subseteq A$ is a unital $C^*$-subalgebra and that $K_A\vert_{B_{\rm sa}} = K_B$, where we implicitly extend the seminorms by $\infty$ outside the finite domains. Suppose $E \colon A \to B$ is a unital completely positive map with $E \vert_{B} = \id_B$. If there is $\varepsilon > 0$ such that for every $a \in \dom (K_A)$
    \begin{align*}
        K_B(E(a)) \leq K_A(a), \quad \Vert a - E(a) \Vert \leq \varepsilon K_A(a),
    \end{align*}
    then
    \begin{align*}
        \Lambda^{\rm dual}((A, K_A), (B, K_B)) \leq \varepsilon. 
    \end{align*}
\end{lemma}
\begin{proof}
    The result will follow from a minor amendment to the proof of \cite[Lemma 2.1]{AguilarLatremoliereRainone}. 
    Since $E$ is a unital completely positive retraction onto $B$, it is a conditional expectation by Tomiyama's theorem.
    The restriction identity gives $K_A(E(a))=K_B(E(a))\leq K_A(a)$, and the approximation hypothesis gives
    \begin{align*}
        \max\{K_A(a),\varepsilon^{-1}\|a-E(a)\|\}=K_A(a).
    \end{align*}
    Thus the cited construction gives a tunnel of extent at most $\varepsilon$, using the $C^*$-algebra $C = A \oplus B$, whose seminorm is
    \begin{align}\label{eq:prop-retraction}
        M(a,b)=\max\{K_A(a),K_B(b),\varepsilon^{-1}\Vert a-b\Vert \}.
    \end{align}
    It remains only to note that this tunnel is Leibniz. Indeed, for self-adjoint $(a,b),(c,d)\in A\oplus B$, the identity $ac-bd=a(c-d)+(a-b)d$ gives 
    \begin{align*}
        \Vert ac - bd \Vert \leq \Vert (a,b) \Vert \Vert c-d \Vert + \Vert (c,d) \Vert \Vert a-b \Vert. 
    \end{align*}
    A straightforward calculation will now verify that $M$ given by \eqref{eq:prop-retraction} is Leibniz. 
\end{proof}

\section{Quantum Codes as Spectral Triples}
\label{sec:quantum-codes-as-spectral-triples}
In this section we show how to construct a unital spectral triple from a stabilizer code. This spectral triple recovers the minimum code distance of the code. For Pauli codes with Hamming weight, it therefore determines correction of all errors up to a prescribed weight.
The starting point for this construction will be a projective error model $(G,\pi)$ and a length function $\ell$ on $G$. As we will be working with this data throughout this section we make the following definition.

\begin{definition}\label{def:projective error model with a weight}
    Let $\calH$ be a Hilbert space. A \emph{projective error model with a weight on $\calH$} is a triple $(G,\pi,\ell)$ where $(G,\pi)$ is a projective error model on $\calH$ and $\ell \colon G \to [0,\infty)$ is a length function.
\end{definition}

Note that since the group $G$ appearing in a projective error model is always assumed to be finite, we have that any length function on $G$ is automatically proper.

There are many sources of length functions to consider. For our purposes we will mostly be interested in the following, which is supposed to resemble a situation where one is aware of the severity or \emph{cost} of a generating set of errors/operations.

\begin{definition}\label{def:cost functions}
    Let $G$ be a finite group and $F \subset G$ be a generating set such that $F = F^{-1}$ and $e \notin F$. A \emph{cost function on $F$} is a function $c \colon F \to (0,\infty)$ such that $c(s^{-1}) = c(s)$ for all $s \in F$.
\end{definition}

The function $c$ is supposed to represent the cost of implementing operations on a physical system, and $c(s) >0$ means every operation costs something, and $c(s) = c(s^{-1})$ means that the cost of undoing the operation is the same as doing the operation to start with. We could artificially extend this definition to the identity $e$ as well, setting $c(e) = 0$, meaning the cost of doing nothing is zero. From this we may define a length function on $G$ as follows.

\begin{definition}\label{def:optimized implementation length}
    Let $G$ be a finite group with a generating set $F$ satisfying $F = F^{-1}$ and $e \notin F$, and let $c \colon F \to (0,\infty)$ be a cost function. The \emph{optimized implementation length} is the length function $\ell_c \colon G \to [0,\infty)$ defined by
    \begin{equation}\label{eqn:optimized implementation length}
        \ell_c(x) = \inf \{ \sum_{i = 1}^N c(s_i) \mid x = s_1 \cdots s_N , s_i \in F \} \quad (x\neq e), \qquad \ell_c(e) = 0.
    \end{equation}
\end{definition}

It is not hard to show that $\ell_c$ is a length function (though symmetry depends on symmetry of $c$, which is going to be immediate for how we have defined a cost function). We also have $\ell_c(x) > 0$ for $x \neq e$ by how we have defined $c$. It is also not hard to see that the infimum may be taken to be a minimum, i.e., a minimizing word exists: Find one word representing $x$ which has finite cost, and note that there are only finitely many words with cost below this. 
There are two immediate uses for this definition:
\begin{itemize}
    \item We may view $c$ as the physical implementation cost which describes how expensive a logical gate is to execute.
    \item We may view $c$ as the error cost of physical errors happening to a code. We would then like to be able to correct all errors below a threshold.
\end{itemize}
These optimized implementation length functions are very satisfying from the physical point of view where one often has an idea of an appropriate cost function for one's system. In particular we may use this idea for the Pauli error model to recover the quantum Hamming weight in the following way.

\begin{example}\label{ex:quantum Hamming weight}
    Let $P_n = (\bbZ_2 \times \bbZ_2)^n$ and $\pi_n \colon P_n \to U((\bbC^2)^{\otimes n})$ be the irreducible projective representation defined in example \ref{ex:projective Pauli error model}, i.e.
    \begin{equation*}
        \pi_n(a_1,b_1,\dots,a_n,b_n) = X^{a_1} Z^{b_1} \otimes \cdots \otimes X^{a_n}Z^{b_n}.
    \end{equation*}
    Let $F \subset P_n$ be the following generating set:
    \begin{equation*}
        F = \bigcup_{i = 1}^{n} \{ (a_1,b_1,\dots,a_n,b_n) \in P_n \mid (a_j,b_j) = (0,0) \text{ for all } j \neq i, \text{ and } (a_i,b_i) \neq (0,0) \}.
    \end{equation*}
    Let $c \colon F \to (0,\infty)$ be the constant function $1$, i.e. $c(s) = 1$ for all $s \in F$. This corresponds to giving the Pauli $X$, the Pauli $Z$, and their product $XZ$ a weight of $1$. The optimized implementation length $\ell_c$ will in this case be the usual quantum Hamming weight, and this gives rise to a projective error model with a weight on $(\bbC^2)^{\otimes n}$.
\end{example}

\begin{construction}\label{const:a logical projective representation}
    Let $(G,\pi)$ be a projective error model on $\calH$ such that $|G| = \dim(\calH)^2$ and let $\calC \subset \calH$ be a stabilizer code. For this construction we will use the notation $S := S(\calC)$ and $L := L(\calC)$. We want to define a projective representation of $\calL_\calC$ on $\calC$. The following construction will yield such a projective representation, though it will depend on a choice of coset representatives for the projective unitary group $PU(\calC)$. To make such a choice, we fix a function $i \colon PU(\calC) \to U(\calC)$ such that $q \circ i = \id_{PU(\calC)}$ and $i(q(I_\calC))=I_\calC$ (recall that $q \colon U(\calC) \to PU(\calC)$ denotes the quotient map).

    Define $\gamma^\prime \colon L \to U(\calC)$ by
    \begin{equation*}
        \gamma^\prime(x) = P_\calC\pi(x)P_\calC.
    \end{equation*}
    By the definition of $L$ and $S$ (see Definition \ref{def:logical group and stabilizer group}) we have that $\gamma^\prime$ is a projective representation of $L$ and that $S = \ker(q \circ \gamma^\prime)$. Hence, the group homomorphism $q \circ \gamma^\prime \colon L \to PU(\calC)$ induces an injective group homomorphism $\tilde{\gamma} \colon \calL_\calC \to PU(\calC)$. Define 
    $\gamma := i \circ \tilde{\gamma} \colon \calL_\calC \to U(\calC)$.
    Since $\tilde{\gamma}$ is an injective group homomorphism and $q \circ \gamma = \tilde{\gamma}$, we have that $\gamma$ is a projectively faithful, projective representation of $\calL_\calC$. We have the following relation between $\gamma$ and $\gamma^\prime$: For any $x \in L$ we have that
    \begin{equation*}
        q(\gamma(xS)) = q(\gamma^\prime(x)).
    \end{equation*}
\end{construction}

The construction outlined in \ref{const:a logical projective representation} motivates the following definition.

\begin{definition}\label{def:the logical algebra}
    Let $(G,\pi)$ be a projective error model on $\calH$ such that $|G| = \dim(\calH)^2$ and let $\calC \subset \calH$ be a stabilizer code. Let $\gamma \colon \calL_\calC \to U(\calC)$ be a projective representation such that
    \begin{equation}\label{eqn:defining equation of gamma}
        q(\gamma(xS(\calC))) = q(P_\calC \pi(x) P_\calC).
    \end{equation}
    The \emph{logical algebra of $\calC$ (with respect to $\gamma$)} is defined to be $A_\calC := \Cred(\calL_\calC,\sigma_\gamma)$.
\end{definition}

Naming $A_\calC$ the logical algebra of $\calC$ is justified by the following result.

\begin{proposition}\label{prop:logical algebra isomorphic to B(C)}
    Let $(G,\pi)$ be a projective error model on $\calH$ such that $|G| = \dim(\calH)^2$ and let $\calC \subset \calH$ be a stabilizer code. With the notation from Definition \ref{def:the logical algebra}, we have that
    \begin{equation*}
        A_\calC \cong B(\calC).
    \end{equation*}
    In particular we have that the isomorphism class of $A_\calC$ is independent of the choice of projective representation $\gamma$.
\end{proposition}
\begin{proof}
    Write $S = S(\calC)$ and $L = L(\calC)$. Since $S$ is normal, we have for any $x \in G$ that $\pi(x) \calC$ is a joint eigenspace for $\pi(S)$. If $x \not \in L$, it is a different joint eigenspace from $\calC$, hence it is orthogonal to $\calC$. Thus $P_\calC \pi(x) P_\calC = 0$ for $x \not \in L$. Irreducibility of $\pi$ gives $\operatorname{span} \pi(G) = B(\calH)$, so
    \begin{align*}
        B(\calC) = \operatorname{span} \{P_\calC \pi(x) P_\calC \mid x \in L \} = \operatorname{span} \gamma(\calL_\calC).
    \end{align*}
    Consequently the integrated representation $\gamma \colon A_\calC \to B(\calC)$ is surjective. By Lemma \ref{lemma:size of pure logical group}, its domain and codomain both have dimension $|\calL_\calC|$, so it is an isomorphism.
\end{proof}

If $(G,\pi,\ell)$ is a projective error model with a weight on $\calH$ such that $|G| = \dim(\calH)^2$, then we want a coherent method of extending the weight $\ell$ from $G$ to all of $B(\calH)$ via the projective representation $\pi$. The spectral triple $(C^*(G,\sigma_\pi),\ell^2(G),D_\ell)$ allows us to do precisely this by considering the resulting Lip-norm $K_\ell \colon C^*(G,\sigma_\pi) \to [0,\infty)$ defined by
\begin{equation*}
    K_\ell(a) = \|[D_\ell,a]\|.
\end{equation*}
Since $\pi$ induces a $*$-isomorphism $C^*(G,\sigma_\pi) \cong B(\calH)$ (the proof of this is analogous to the proof of Proposition \ref{prop:logical algebra isomorphic to B(C)}), this gives us a coherent definition of weight of an element of $B(\calH)$. In particular we have that
\begin{equation}\label{eqn:weight of group element}
    K_{\ell}(\lambda^{\sigma_\pi}_s)
    = \|[D_\ell,\lambda^{\sigma_\pi}_s]\|
    = \sup\{ \vert \ell(st) - \ell(t) \vert \mid t \in G \}
    = \ell(s),
\end{equation}
which shows that $K_\ell$ really is an extension of $\ell$ from $G$ to all of $\Cred(G,\sigma_\pi) \cong B(\calH)$.

Note that the pair $(C^*(G,\sigma_\pi),K_\ell)$ is always a compact quantum metric space whenever $G$ is finite. We will now outline how to construct a compact quantum metric space associated to a stabilizer code whenever $|G| = \dim(\calH)^2$: Suppose that $\calC \subset \calH$ is a stabilizer code. Since $(G,\pi,\ell)$ is a projective error model with a weight on $\calH$, we get an induced length function $\ell_\calC \colon \calL_\calC \to [0,\infty)$ by setting
\begin{equation}\label{eqn:length function on logical algebra}
    \ell_\calC(xS(\calC)) = \min\{ \ell(xs) \mid s \in S(\calC) \}.
\end{equation}
We will refer to this as the \emph{quotient length function} or \emph{quotient weight}.

\begin{definition}\label{def:the logical spectral triple}
    Let $(G,\pi,\ell)$ be a projective error model with a weight on $\calH$ such that $|G| = \dim(\calH)^2$ and let $\calC \subset \calH$ be a stabilizer code. The \emph{logical spectral triple of $\calC$} is the unital spectral triple $(A_\calC, \ell^2(\calL_\calC), D_{\calC})$, where $D_\calC := D_{\ell_\calC}$. We will denote the resulting Leibniz quantum compact metric space by $(A_\calC,K_\calC)$, where  $K_\calC :=K_{\ell_\calC}\vert_{(A_\calC)_{\rm sa}}$.
\end{definition}

\begin{remark}
    So far we have suppressed the projective representation $\gamma$ from the notation of the logical spectral triple $(A_\calC, \ell^2(\calL_\calC), D_{\calC})$, Proposition \ref{prop:logical algebra isomorphic to B(C)} justifies that we may drop $\gamma$ from the notation of the logical algebra $A_\calC$. In the following definition we will make the dependence of $\gamma$ explicit, but Theorem \ref{thm:minimum code distance is minimum of basis logicals} proves that  
    key properties of the code 
    we care about, namely its error correcting abilities and the isomorphism class of $A_\calC$, are independent of the projective representation $\gamma$.
\end{remark}

\begin{definition}\label{def:code distance}
    Let $(G,\pi,\ell)$ be a projective error model with a weight on $\calH$ such that $|G| = \dim(\calH)^2$. If $\calC \subset \calH$ is a stabilizer code, we define its \emph{minimum code distance (with respect to $(G,\pi,\ell)$)} to be the number
    \begin{equation*}
        \Delta_\calC = \inf\{ K_{\ell_\calC}(a) \mid a \in A_\calC = C^*(\calL_\calC,\sigma_\gamma), \ \tau_{\sigma_\gamma}(a^*a) = 1, \text{ and } \tau_{\sigma_\gamma}(a) = 0 \}.
    \end{equation*}
    We use the convention $\inf\varnothing=+\infty$. In particular, $\Delta_\calC=+\infty$ when $\dim\calC=1$.
\end{definition}

This definition is theoretically satisfying as it looks at all the possible (normalized) logical operators that can act on the code $\calC$ in a non-trivial manner and insists that the minimum code distance is the minimum weight of all of these operators. This is completely analogous to how the minimum code distance can be computed for classical linear codes. The infimum appearing in Definition \ref{def:code distance} makes it so that $\Delta_\calC$ is a priori difficult to compute (it will turn out that computing $\Delta_\calC$ is always difficult, but not as difficult as this definition would make you think). Luckily the following result tells us that we only have to compute the weight of the elements coming from the group $\calL_\calC$.

\begin{theorem}\label{thm:minimum code distance is minimum of basis logicals}
    Let $(G,\pi,\ell)$ be a projective error model with a weight on $\calH$ such that $|G| = \dim(\calH)^2$. If $\calC \subset \calH$ is a stabilizer code, then
    \begin{equation*}
        \Delta_\calC = \min\{ \ell_\calC(x) \mid x \in \calL_\calC \setminus \{e\} \}.
    \end{equation*}
    Here the minimum over the empty set is understood to be $+\infty$. 
    In particular, $\Delta_\calC$ is independent of the choice of projective representation $\gamma \colon \calL_\calC \to U(\calC)$ satisfying \eqref{eqn:defining equation of gamma}.
\end{theorem}
\begin{proof}
    If $\dim\calC=1$, both sides are $+\infty$. Assume therefore in the sequel that $\dim\calC>1$.
    Let $a \in A_\calC = C^*(\calL_\calC,\sigma_\gamma)$ be any element satisfying $\tau_{\sigma_\gamma}(a^*a) = 1$ and $\tau_{\sigma_\gamma}(a) = 0$. We may write
    \begin{equation*}
        a = \sum_{x \in \calL_\calC} a_x \lambda^{\sigma_\gamma}_x.
    \end{equation*}
    In this case we have that
    \begin{equation*}
        \sum_{x \in \calL_\calC} |a_x|^2 = \tau_{\sigma_\gamma}(a^*a) = 1 \quad \text{and} \quad a_{e} = \tau_{\sigma_\gamma}(a) = 0.
    \end{equation*}
    For ease of notation we write
    \begin{equation*}
        \Delta := \min\{ \ell_\calC(x) \mid x \in \calL_\calC \setminus \{e\} \}.
    \end{equation*}
    By the observation made in \eqref{eqn:weight of group element} specialized to our current situation, we have that $K_{\ell_\calC}(\lambda^{\sigma_\gamma}_x) = \ell_\calC(x)$. This ensures that $\Delta_\calC \leq \Delta$. With this we compute the following:
    \begin{align*}
        \Delta_\calC^2
        \leq \Delta^2
        = \sum_{x \in \calL_\calC} |a_x|^2\Delta^2 
        \leq \sum_{x \in \calL_\calC} |a_x|^2\ell_\calC(x)^2 
        = \|[D_{\ell_\calC},a]\delta_{e}\|^2 
        \leq K_{\ell_\calC}(a)^2.
    \end{align*}
    The second inequality follows because $a_{e} = 0$. As $a$ was arbitrarily chosen, we have that $\Delta_\calC = \Delta$.
\end{proof}

In the special case where we consider the Pauli error model and the optimized implementation length as in example \ref{ex:quantum Hamming weight}, then Theorem \ref{thm:minimum code distance is minimum of basis logicals} tells us that we recover the standard notion of minimum code distance used for Pauli stabilizer codes with $\dim \calC > 1$ coming from the quantum Hamming weight. We also recover the following error correction result.

\begin{proposition}\label{prop:error correction in terms of distance}
    Let $(G,\pi,\ell)$ be a projective error model with a weight on $\calH$ such that $|G| = \dim(\calH)^2$ and let $\calC \subset \calH$ be a stabilizer code. Then $g \in G$ satisfies the Knill--Laflamme conditions, i.e.
    \begin{equation*}
        P_\calC \pi(g) P_\calC \in \bbC P_\calC
    \end{equation*}
    if $K_\ell(\lambda_g^{\sigma_\pi}) < \Delta_\calC$.
\end{proposition}
\begin{proof}
    By \eqref{eqn:weight of group element} we have that $K_\ell(\lambda_g^{\sigma_\pi}) = \ell(g)$. Suppose that $\ell(g) < \Delta_\calC$. There are two cases to consider, namely $g \in L(\calC)$ and $g \in G \setminus L(\calC)$. In the first case we have that
    \begin{equation*}
        \ell_\calC(gS(\calC)) \leq \ell(g) < \Delta_\calC
    \end{equation*}
    which implies that $g \in S(\calC)$ by Theorem \ref{thm:minimum code distance is minimum of basis logicals}. In particular we have that $P_\calC \pi(g) P_\calC \in \bbT P_\calC$ in this case.

    In the other case, if $g \in G \setminus L(\calC)$, we have that $P_\calC \pi(g) P_\calC = 0$ by \cite[Theorem 6.3]{Eidesen2026}.
\end{proof}

We end this section by highlighting the following definition of an $[[n,k,d]]_p$-Pauli stabilizer code which we will work with for the remainder of this paper. Setting $p = 2$ in the following recovers what we studied in examples \ref{ex:projective Pauli error model} and \ref{ex:quantum Hamming weight}.

\begin{definition}\label{def:[[n,k,d]]-codes}
    Let $p \geq 2$ be a prime, $n\geq 0$ and write $\bbF_p := \bbZ_p$. Set $P_{p,n} = \bbF_p^n \oplus \bbF_p^n$, and define $\pi_{p,n} \colon P_{p,n} \to U((\bbC^p)^{\otimes n})$ by
    \begin{equation*}
        \pi_{p,n}(x,z) = X_p^{x_1}Z_p^{z_1} \otimes \cdots \otimes X_p^{x_n}Z_p^{z_n}
    \end{equation*}
    where
    \begin{align*}
        X_p |j\rangle = |j+1 \rangle, \qquad Z_p |j\rangle = \zeta^j |j\rangle, \qquad \zeta = e^{2\pi i/p}.
    \end{align*}
    Define the Hamming weight $\ell_{p,n} \colon P_{p,n} \to [0,\infty)$ by
    \begin{equation*}
        \ell_{p,n}(x,z) = |\{ i \mid (x_i,z_i) \neq (0,0) \}|.
    \end{equation*}
    The projective error model with a weight on $(\bbC^p)^{\otimes n}$, $(P_{p,n}, \pi_{p,n}, \ell_{p,n})$, is called the \emph{$p$-qudit Pauli error model}. We will refer to elements of $P_{p,n}$ as \emph{labels}.
    
    For integers $1 \leq k \leq n$ and $0 \leq d$, an \emph{$[[n,k,d]]_p$-Pauli stabilizer code} is a $p^k$-dimensional stabilizer code $\calC \subset (\bbC^p)^{\otimes n}$ (with respect to the Pauli error model $(P_{p,n},\pi_{p,n})$) such that $\Delta_\calC = d$.
\end{definition}

\begin{remark}\label{rem:Pauli stabilizer codes always yield [[n,k,d]]-codes}
    In the case where $p = 2$, it is well-known throughout the literature (see for example \cite{NielsenChuang2010} for a textbook reference) that every Pauli stabilizer code $\calC$ has dimension $2^k$ for some integer $0 \leq k \leq n$, and computing the number $k$ is relatively straightforward given the subgroup $S \subset P_n$ that gives rise to the stabilizer code $\calC$. If $k\geq 1$, Theorem \ref{thm:minimum code distance is minimum of basis logicals} also shows that $\Delta_\calC$ is an integer, hence every such Pauli stabilizer code is an $[[n,k,d]]_2$-Pauli stabilizer code for some $k$ and $d$. The hard part is in general to compute the number $\Delta_\calC$, this is in fact known to be NP-hard, cf. \cite{KapshikarKundu2023}.
\end{remark}

\section{Quantum distances between \texorpdfstring{$[[n,k,d]]_p$}{[[n,k,d]]}-codes}\label{sec:quantum-distances-between-codes}
In this section we compare $[[n,k,d]]_p$-codes $\calC$ through the quantum geometry of their corresponding Leibniz quantum compact metric spaces $(A_\calC,K_\calC)$ using the dual propinquity from Section \ref{subsec:propinquity}. Everything in this section will therefore be done relative to the projective error model with a weight $(P_{p,n},\pi_{p,n},\ell_{p,n})$, and we will not mention this explicitly again in this section. We start by recalling some useful facts about stabilizer codes in this setting, as well as establishing some notation that will be useful.

\begin{lemma}
    Let $\calC$ be a stabilizer code. Define the commutator form $\omega \colon P_{p,n} \times P_{p,n} \to \bbF_p$ by
    \begin{equation*}
        \omega ((x,z) ,(x', z'))  = z \cdot x' - z' \cdot x
    \end{equation*}
    where the dot product is taken in the vector space $\bbF_p^n$. Then, setting $S := S(\calC)$, we have that
    \begin{equation*}
        L(\calC) = S^\perp = \{ v \in P_{p,n} \mid \omega(v,s) = 0 \text{ for all } s \in S \}.
    \end{equation*}
\end{lemma}
\begin{proof}
     Note that $\sigma_{\pi_{p,n}}((x,z),(x',z')) = \zeta^{z \cdot x'}$ where $\zeta = e^{2\pi i /p}$. By point (1) in Proposition \ref{prop:collection of stabilizer code results} we have that
     \begin{align*}
         L(\calC) 
         & = \{ (x,z) \in P_{p,n} \mid \lambda_\calC(s_x,s_z) = \zeta^{z \cdot s_x}\zeta^{- s_z \cdot x}\lambda_\calC(s_x,s_z) \text{ for all }(s_x,s_z) \in S \} \\
         & = \{ v \in P_{p,n} \mid \omega(v,s) = 0 \text{ for all } s \in S \}. \qedhere
     \end{align*}
\end{proof}

By point (4) in Proposition \ref{prop:collection of stabilizer code results} we also have that $\dim(\calC) = p^k$ where $k = n - \dim_{\bbF_p}(S)$, hence
\begin{equation}\label{eqn:boundary logical algebra}
    \calL_{\calC} = S^\perp/S\cong\bbF_p^{2k}, \qquad A_{\calC} \cong M_{p^k}(\bbC).
\end{equation}
Consider the quotient length $\ell_{\calC}$ of \eqref{eqn:length function on logical algebra}. By Theorem \ref{thm:minimum code distance is minimum of basis logicals}, $\Delta_{\calC}$ is the least weight of a label in $S^\perp\setminus S$.

\subsection{A topology on stabilizer codes}

\begin{definition}\label{def:code-distance-and-boundary}
    Consider $[[n,k,d]]_p$-codes $\calC$ and $\calC'$ and their Leibniz quantum compact metric spaces $(A_\calC,K_\calC)$ and $(A_{\calC'},K_{\calC'})$ of Definition \ref{def:the logical spectral triple}. Define a pseudometric $d_p$ between two such codes by
    \begin{align*}
        d_p (\calC, \calC') := \Lambda^{\rm dual} ((A_\calC, K_\calC), (A_{\calC'}, K_{\calC'}))
    \end{align*}
    using \eqref{eq:def-dual-propinquity}. Let $\calQ_p$ denote the set of pairs $(A_\calC,K_\calC)$ where $\calC$ is an $[[n,k,d]]_p$-code, up to full quantum isometry. That is, $(A_\calC,K_\calC) \sim (A_{\calC^\prime},K_{\calC^\prime})$ if and only if there exists a $*$-isomorphism $\theta \colon A_\calC \to A_{\calC^\prime}$ such that $K_{\calC^\prime}(\theta(a)) = K_\calC(a)$ for all $a \in A_{\calC, {\rm sa}}$. Its \emph{metric boundary} is $\partial \calQ_p = \overline{\calQ_p} \setminus \calQ_p$, where the closure is taken in the space of Leibniz quantum compact metric spaces. 
\end{definition}

Whenever calculating propinquity we restrict the seminorm  to the self-adjoint part. Changing the projective representation $\gamma$ which is implicit in the definition of $(A_\calC, K_\calC)$  gives a diagonal unitary equivalence between the regular representations which commutes with the Dirac operator, so the assignment is independent of this choice. The distance $d_p$ is a pseudometric on codes and a metric on $\calQ_p$. At distance zero the logical algebras are isomorphic, so the dimensions of the corresponding codes agree by Proposition \ref{prop:logical algebra isomorphic to B(C)}. The number $n$ of tensor factors need not agree, as replacing $\calC$ by $\calC \otimes \bbC |0\rangle$ leaves the quotient length $\ell_\calC$ unchanged. 

We record the following lemma which is particularly useful as the codes under consideration arise from finite abelian groups. 

\begin{lemma}\label{lemma:fourier-truncation}
    Let $G$ be a finite abelian group with length $\ell$ and $2$-cocycle $\sigma$. For a subgroup $H \leq G$, let $E_H \colon C^*_r(G,\sigma) \to C^*_r(G, \sigma)$ be the map which discards the Fourier coefficients outside $H$. It is a unital completely positive retraction, which is moreover contractive for the commutator seminorm. If $\ell (g) \geq R >0$ for $g \not \in H$, we have
    \begin{align}
        \Vert a - E_H (a) \Vert \leq \frac{\sqrt{|G \setminus H|}}{R} K_\ell (a). 
    \end{align}
\end{lemma}

\begin{proof}
    The set $H^\circ = \{\chi \in \widehat{G} \mid \chi(h) = 1 \text{ for all } h \in H \}$ is the annihilator of $H$, and coincides with the characters factoring through $G/H$. Each $\chi \in \widehat{G}$ defines a $*$-automorphism on $C^*_r(G, \sigma)$ through
    \begin{align*}
        \alpha_\chi (\lambda_g^\sigma) = \chi (g) \lambda_g^\sigma. 
    \end{align*}
    This is spatially implemented on $\ell^2(G)$ through the diagonal unitary $U_\chi \delta_g = \chi(g) \delta_g$. We observe that $U_\chi$ commutes with $D_\ell$, and so $K_\ell(\alpha_\chi(a)) = K_\ell(a)$ for all $a \in C^*_r (G, \sigma)$. 
    Now with $E_H := |H^\circ |^{-1} \sum_{\chi \in H^\circ} \alpha_\chi$ we have
    \begin{align*}
        E_H (\lambda^\sigma_g) = \begin{cases}
            \lambda^\sigma_g & g\in H \\
            0 & g \not\in H.
        \end{cases}
    \end{align*}
    Thus $E_H$ is a unital completely positive retraction onto the subgroup algebra, and
    \begin{align*}
        K_\ell(E_H(a)) \leq \frac{1}{|H^\circ|} \sum_{\chi \in H^\circ} K_\ell(\alpha_\chi (a)) = K_\ell(a). 
    \end{align*}
    Using the standard argument for $a = \sum_g a_g \lambda_g^\sigma$ we get
    \begin{align}\label{eq:lower-2-bound}
        \sum_{g \in G} \vert a_g \vert^2 \ell(g)^2 = \Vert [D_\ell, a] \delta_e \Vert^2 \leq K_\ell(a)^2,
    \end{align}
    so for $R >0$ satisfying $\ell(g) \geq R$ for $g \not \in H$, we use Cauchy-Schwarz to obtain
    \begin{equation*}
        \Vert a - E_H (a) \Vert 
        \leq \sum_{g \not \in H} \vert a_g \vert 
        \leq \sum_{g \not \in H} \vert a_g \vert \frac{\ell(g)}{R} 
        \leq \frac{\sqrt{|G \setminus H|}}{R} \bigg( \sum_{g\not \in H} \vert a_g \vert^2 \ell(g) ^2 \bigg)^{1/2} 
        \leq \frac{\sqrt{|G \setminus H |}}{R} K_\ell(a). \qedhere
    \end{equation*}
\end{proof}

\subsubsection{Examples of convergence}
\begin{example}\label{example:repetition-code}
    In $(\bbC^p)^{\otimes N}$, that is on $N \geq 1$ qudits, let $\mathcal{R}_N := \operatorname{span} \{ | j \rangle^{\otimes N} \mid j \in \bbF_p\}$, and let $Z_{p,i}$ act on the $i$th factor. For $1 \leq i < N$, the operators $T_i = Z_{p,i}Z_{p,i+1}^{-1}$ satisfy
    \begin{align*}
        T_i |j_1 \ldots j_N \rangle = \zeta^{j_i - j_{i+1}} |j_1 \ldots j_N \rangle
    \end{align*}
    for $\zeta = e^{\frac{2\pi i}{p}}$. 
    Thus the joint $+1$ eigenspace necessitates $j_1 = j_2 = \cdots = j_N$, meaning it agrees with $\mathcal R_N$. The corresponding stabilizer label space is $S_N = \{ (0,z) \mid \sum_i z_i = 0 \}$. 

    The commutation relation $T_i \pi_{p,N}(x,z) = \zeta^{x_i - x_{i+1}} \pi_{p,N}(x,z) T_i$ implies $S_N^\perp = \{(a \mathbf 1, z) \mid a \in \bbF_p, z \in \bbF_p^N \}$ for $\mathbf 1 = (1, \ldots , 1)$. The map
    \begin{align*}
        q \colon S_N^\perp \to \bbF_p^2, \qquad q(a \mathbf 1, z) = (a , \sum_i z_i)
    \end{align*}
    has kernel $S_N$ and identifies $S_N^\perp /S_N$ with $\bbF_p^2$. 
    If $b = \sum_i z_i$, the restricted action is
    \begin{align*}
        \pi_{p,N}(a \mathbf 1, z) |j\rangle^{\otimes N} = \zeta^{bj} |j+a \rangle^{\otimes N},
    \end{align*}
    and so we may choose logical generators $\overline{X} = X_p^{\otimes N}\vert_{\mathcal R_N}$ and $\overline{Z} = Z_{p,1}\vert_{\mathcal R_N}$. The corresponding Weyl operators $u_{a,b} = \overline{X}^a \overline{Z}^b$ have cocycle $\sigma((a,b), (a',b')) = \zeta^{ba'}$. 

    The code has dimension $p$  and encodes only one logical qudit. For $a \neq 0$, every representative of $(a,b)$ has a nonzero shift in every tensor leg. For $a= 0$ and $b \neq 0$, the representative $Z_{p,1}^b$ attains the least possible weight. Thus
    \begin{align*}
        \ell_N (0,b) &= 1 \qquad \text{for $b \neq 0$} \\ 
        \ell_N(a,b) &= N \qquad \text{for $a \neq 0$}
    \end{align*}
    In particular, the code parameters are $[[N,1,1]]_p$. 

Now put $J = \{0\} \times \bbF_p$, and note that its subgroup algebra $B = C^*_r(J, \sigma|_J) \cong C(\bbZ_p)$. Moreover, $B$ carries the seminorm $K_J$ coming from the length function which assigns the value $1$ to every element of $J \setminus \{0\}$. 
Indeed, decompose the regular space as $\ell^2(\bbF_p^2 ) = \oplus_{c \in \bbF_p} \calH_c$ for $\calH_c = \operatorname{span} \{\delta_{(c,d)} \mid d \in \bbF_p \}$. Then
\begin{align*}
    \lambda_{(0,b)}^\sigma\delta_{(c,d)} =\zeta^{bc}\delta_{(c,d+b)},\qquad D_{\ell_N}|_{\mathcal H_c} =
    \begin{cases}
 I-P_{\delta_{(0,0)}},&c=0,\\ 
 NI,&c\ne0.
 \end{cases}
\end{align*}
    Each $\calH_c$ is $B$-invariant. For $c=0$ we get the regular commutator on $J$, while $c\neq 0$ gives zero as the Dirac is scalar on these blocks. This verifies $K_{\mathcal R_N}|_{B_{\rm sa}}=K_J$ and we may employ Lemma \ref{lemma:fourier-truncation} and Lemma \ref{lemma:cp-approximation-map} to deduce that
    \begin{align*}
         \Lambda^{\rm dual}\big((A_{\mathcal R_N},K_{\mathcal R_N}),   (C(\bbZ_p),K_J)\big) \le\frac{\sqrt{p^2-p}}{N}.
    \end{align*}
    Thus the conditional expectations onto $B$ approximate the identity uniformly on the respective Lip-unit balls, and the associated quantum compact metric spaces converge to $(C(\bbZ_p), K_J)$.
\end{example}

To construct limits with matrix algebra factors, we first realize $M_p(\bbC)$ with a quotient length that is constant on the nonzero labels. Iterating this construction will give the scales used in the AF algebra approximations below.

\begin{lemma}\label{lemma:513-protector-code}
    For every prime $p$ there exists a $[[5,1,3]]_p$-stabilizer code $\calF$ whose quotient Hamming length is three on every nonzero element of $\calL_{\calF}=L(\calF)/S(\calF)$. 
    For each $j\ge1$, there exists a code $\calF^{(j)}$ with parameters $[[5^j,1,3^j]]_p$ and which satisfies
    \begin{align*}
        \ell_{\calF^{(j)}}(x)=3^j\qquad \text{for every }x\in\calL_{\calF^{(j)}}\setminus\{0\}.
    \end{align*}
\end{lemma}

\begin{proof}
    Rains's stabilizer construction \cite[Theorem 13]{RainsNBQC99} yields a $p$-dimensional code on 5 $p$-dimensional sites (where \emph{site} refers to one of the tensor factors appearing in $(\bbC^p)^{\otimes 5}$) with minimum distance $3$, that is, a $[[5,1,3]]_p$-code in our notation. 
    Put $P = P_\calF$, $S=S(\calF)$, and $\pi = \pi_{p,5}$.
    For any two sites $T$, let $\calA_T$ be the full matrix algebra on those sites, acting as the identity on $T^c$. As $\calA_T$ can be written as the span of elements of the form $\pi(g)$ where $\ell_{p,5}(g) < 3 = \Delta_\calF$, we have $P \calA_T P = \bbC P$ by Proposition \ref{prop:error correction in terms of distance}. 
    Each $\pi(s)$ normalizes $\calA_T$, so
    \begin{align*}
        \calM = C^*(\calA_T, \pi(S)) = \operatorname{span} \{ a \pi(s) \mid a \in \calA_T, s\in S \}. 
    \end{align*}
    Moreover, $P \in C^*(\pi(S))$, and each $\pi(s)$ is scalar on $\calF$. Hence $P \calM P = \bbC P$, and $P$ is a minimal projection of $\calM$. The structure theorem for finite-dimensional $C^*$-algebras gives
    \begin{align*}
        P \calM' P = B(\calF), \qquad \calM' = \calA_{T^c} \cap \pi(S)'.
    \end{align*}
    The latter algebra is spanned by the Pauli operators supported on $T^c$ that commute with $\pi(S)$. Their compressions therefore span $B(\calF)$. 
    Since the operators corresponding to distinct classes in $\calL_\calF$ are linearly independent, every class has a representative supported on $T^c$. Thus every nonzero class has length at most three, and the minimum distance gives the reverse inequality.

    We also need to understand how the quotient lengths behave under concatenation. Let $\calE \subseteq (\bbC^p)^{\otimes n}$ be a $p$-dimensional stabilizer code such that $\ell_\calE (y) = w$ for every nonzero $y \in \calL_\calE$. Choose an isometry $V \colon \bbC^p \to \calE$ which intertwines $X_p$ and $Z_p$ with restrictions of physical Pauli operators, allowing scalar phases. This identifies $\calL_\calE$ with $\bbF_p^2$. The standard stabilizer concatenation construction sends a stabilizer code $\calC \subseteq (\bbC^p)^{\otimes N}$ to $\calC[\calE] = V^{\otimes N}\calC$, and identifies their logical label groups by replacing the Pauli operators on each site with the corresponding operators on $\calE$, see \cite[Proposition 1]{WangZengGrasslSanders}. Under this identification, writing $u_i$ for the $i$th coordinate pair of $u \in \bbF_p^N \oplus \bbF_p^N$, we obtain
    \begin{align*}
        \ell_{\calC[\calE]}(x) = \min_{\substack{u\in L(\calC)\\u+S(\calC)=x}} \sum_{i=1}^N \ell_{\calE} (u_i) = w \ell_\calC (x). 
    \end{align*}
    Indeed, a Pauli operator either preserves a stabilizer code or has zero compression to it. Thus a Pauli operator preserving $\calC[\calE] \subseteq \calE^{\otimes N}$ must preserve $\calE$ in each block, because otherwise its compression to $\calE^{\otimes N}$ and hence to $\calC[\calE]$ would vanish. Its blockwise restrictions determine a representative $u$ of the corresponding class $x$ for $\calC$, and its weight is at least $\sum_i \ell_\calE(u_i)$. Conversely, for any representative $u$ of this prescribed class, choose a minimum weight representative of each inner class $u_i$. Their tensor product acts on $\calE^{\otimes N}$ as $V^{\otimes N} \pi_{p,N}(u) (V^*)^{\otimes N}$ up to a scalar, and hence represents $x$ for $\calC[\calE]$, with weight exactly $\sum_i \ell_\calE(u_i)$. Minimizing then proves the first equality. The assumed form of $\ell_\calE$ proves the second. 

    Set $\calF^{(1)} = \calF$ and $\calF^{(j+1)} = \calF[\calF^{(j)}]$. Inductively, these codes have dimension $p$, occupy $5^j$ sites, and have quotient length $3^j$ on every nonzero class. Their minimum distance is therefore $3^j$, as required. 
\end{proof}

\begin{proposition}\label{prop:concatentation-qudit-collapse}
    Let $\calC_0$ be a fixed $[[n_0, k_0, d_0 ]]_p$ stabilizer code with $k_0 \geq 1$, and put $\calC_m = \calC_0 \otimes \calF^{(m)}$ for $m\geq 1$. Then $\calC_m$ has parameters $[[n_0 + 5^m , k_0 +1 , \min\{d_0, 3^m \} ]]_p$, and 
    \begin{align}
        d_p (\calC_m , \calC_0) \leq (p^2 -1) 3^{-m}. 
    \end{align}
    Moreover, the pairs associated to $\calF^{(m)}$ converge to $(\bbC, 0)$. 
\end{proposition}

\begin{proof}
    Set $w = 3^m$. 
    A logical class of the product code $\calC_m=\calC_0\otimes\calF^{(m)}$ is specified by a pair $(g,u)\in\calL_{\calC_0}\oplus\calL_{\calF^{(m)}}$. Its physical representatives consist of a representative of $g$ in the first block and a representative of $u$ in the second. These representatives can be chosen independently, and their Hamming weights add. Hence choosing a minimum-weight representative in each block minimizes the total weight, so \eqref{eqn:length function on logical algebra} gives $\ell_{\calC_m}(g,u) =\ell_{\calC_0}(g)+w\mathbf1_{\{u\neq 0\}}$. 
    The code dimensions multiply, and Theorem \ref{thm:minimum code distance is minimum of basis logicals} gives the asserted parameters.

    Choose the tensor product of the logical projective representations from Construction~\ref{const:a logical projective representation}.
    Writing $A_0=A_{\calC_0}$, $D_0=D_{\ell_{\calC_0}}$, and $\sigma$ for the cocycle of the second factor, the unitary $\delta_{(g,u)}\mapsto\delta_g\otimes\delta_u$ identifies the logical spectral triple with
    \begin{align*}
        A_{\calC_m}=A_0\otimes B,\qquad B=\Cred(\bbF_p^2,\sigma),\qquad D=D_0\otimes I+I\otimes w(I-P_{\bbC\delta_0}).
    \end{align*}
    In particular, we have
    \begin{align*}
        [D, b\otimes I] = [D_0, b] \otimes I,
    \end{align*}
    and so the seminorm restricts exactly to $K_{\calC_0}$ on $A_0\otimes I$.

    Let $\tau_\sigma$ be the canonical trace on $B$ and consider the conditional expectation
    \begin{align*}
        E(a) = E_0 (a) \otimes I, \quad E_0 = \id \otimes \tau_\sigma .
    \end{align*}
    For $v \in \bbF_p^2$, let $V_v \xi = \xi \otimes \delta_v$. Then in the left regular representation $E_0(a) = V_0^* a V_0$, while $DV_0 = V_0 D_0$. Thus
    \begin{align*}
        [D_0, E_0(a) ] = V_0^* [D, a] V_0, \quad K_{\calC_m} (E(a)) = K_{\calC_0}(E_0 (a)) \leq K_{\calC_m} (a)
    \end{align*}
    for any $a = a^* \in A_{\calC_m}$. 

    Now for $a=a^*$, write $a = \sum_{u \in \bbF_p^2} a_u \otimes \lambda_u^\sigma$, and note that for $u \neq 0$ two matrix corners of the commutator give
    \begin{align*}
        V_u^* [D,a]V_0 &= [D_0, a_u] + wa_u, \\
        \overline{\sigma(u,-u)} V_0^* [D,a]V_{-u} &= [D_0, a_u] - wa_u. 
    \end{align*}
    Both corners have norm at most $K_{\calC_m}(a)$. Subtracting thus yields $w \Vert a_u \Vert \leq K_{\calC_m}(a)$. As $E(a) = a_0 \otimes I$, it follows that
    \begin{align*}
        \Vert a - E(a) \Vert \leq \sum_{u \neq 0} \Vert a_u \Vert \leq \frac{p^2-1}{w} K_{\calC_m}(a). 
    \end{align*}
    As the data satisfy the hypotheses of Lemma \ref{lemma:cp-approximation-map} we obtain
    \begin{align*}
        d_p (\calC_m, \calC_0 ) \leq (p^2 -1) 3^{-m}. 
    \end{align*}

    For the second factor alone, Lemma \ref{lemma:fourier-truncation} with $H=\{0\}$ and $R=w$ gives, for $b=b^* \in B$,
    \begin{align*}
        \|b-\tau_\sigma(b)I\| \le\frac{\sqrt{p^2-1}}{w}K_{\calF^{(m)}}(b).
    \end{align*}
    Applying Lemma \ref{lemma:cp-approximation-map} using $(\bbC, 0)$, yields
    \begin{align*}
        \Lambda^{\rm dual} ((A_{\calF^{(m)}}, K_{\calF^{(m)}}), (\bbC, 0)) \leq \sqrt{p^2 -1} \cdot 3^{-m} \longrightarrow 0. 
    \end{align*}
    This finishes the proof. 
\end{proof}

The additional  logical qudit remains physically encoded throughout this sequence. Its disappearance from the limit means that it is invisible to observables of bounded logical cost. In particular, taking $k_0=1$ gives convergence from two logical qudits to one.

\subsection{The metric boundary}
For a code $\calC$,  Lemma \ref{lemma:fourier-truncation} and its proof give an isometric action of $\widehat{\calL_\calC}$ on $A_\calC \cong M_{p^k} (\bbC)$. 
The nondegenerate Pauli commutator pairing identifies this dual action with the conjugation by the logical unitaries, $a\mapsto\lambda_g^\sigma a(\lambda_g^\sigma)^*$. 
The action is ergodic. Indeed, applying Lemma \ref{lemma:fourier-truncation} using $H = \{e \}$, the average over all characters sends $a = \sum_{g \in \calL_\calC} a_g \lambda_g^\sigma$ to $a_e 1$. As every fixed element under the action will equal its average, we deduce that such elements are scalar. 
We use this to characterize the possible limits of codes in the following theorem. 
\begin{theorem}\label{thm:limits-are-AF}
Suppose that $(A_j,K_j)$ converges to $(A,K)$ in dual propinquity, where $A_j=M_{p^{k_j}}(\bbC)$ and $K_j$ is a closed Leibniz Lip-norm invariant under the standard logical Pauli adjoint action. Then
 \begin{equation}\label{eq:boundary-algebra-classification}
  A\cong C((\bbZ_p)^r)\otimes B_{p,s}, \qquad r,s\in\{0,1,2,\ldots,\infty\},
 \end{equation}
 where
 \begin{equation*}
  B_{p,s}=
  \begin{cases}
   M_{p^s}(\bbC),&s<\infty,\\
   \displaystyle\bigotimes_{j=1}^\infty M_p(\bbC),&s=\infty.
  \end{cases}
 \end{equation*}
 The notation $(\bbZ_p)^\infty$ denotes a countably infinite product equipped with the product topology, and $B_{p,0}=\bbC$. In particular, every algebra underlying a point of
 $\overline{\calQ_p}$ is AF.
\end{theorem}

We prove the theorem in two steps. Every algebra $A_i$ appearing in a sequence in Theorem \ref{thm:limits-are-AF} comes equipped with an action of some group $G_i$ through Lemma \ref{lemma:fourier-truncation} and its proof. The first of the following lemmas tells us that the limit $A$ also has an ergodic action by a compact metrizable group. The second lets us know that the existence of such an action forces the limit $C^*$-algebra to have a very particular form. 

\begin{lemma}\label{lemma:limit-group-action}
    Let $p$ be a fixed prime. Let $(A_j, K_j)$ be Leibniz quantum compact metric spaces with $A_j = M_{p^{k_j}} (\bbC)$, $k_j \geq 0$, and suppose $K_j$ is invariant under the logical Pauli conjugation action. If $\Lambda^{\rm dual}((A_j,K_j), (A,K))\to 0$, then there is a compact metrizable abelian group $H$ with $h^p = e$ for every $h \in H$, which acts strongly continuously and ergodically on $A$ by full quantum isometries. 
\end{lemma}

\begin{proof}
    Set $H_j \cong (\bbZ_p)^{2k_j}$ with action $\alpha^{(j)}$ on $A_j $, where $A_j$ can also be realized as a twisted group $C^*$-algebra. Set
    \begin{align*}
        X_j = (S(A_j) , \mk_{K_j}), \qquad X = (S(A), \mk_K). 
    \end{align*}
    As observed in Section \ref{subsec:propinquity}, a tunnel of extent $\varepsilon$ implies the Hausdorff distance between the endpoint state spaces is at most $\varepsilon$, and so we deduce $X_j \to X$ in Gromov--Hausdorff distance. 

    Note that we have an action of $H_j$ on $X_j$ through $g \cdot \phi := \phi \circ \alpha_{g^{-1}}^{(j)}$. This action can be seen to be isometric, and so $H_j \leq \operatorname{Isom}(X_j)$. The associated uniform metric is
    \begin{align*}
        d_j (g,h) = \sup \{ \Vert \alpha_g^{(j)} (a) - \alpha_h^{(j)}(a) \Vert \mid a=a^*, K_j(a) \leq 1 \}.
    \end{align*}
    We recall that the norm of a self-adjoint element $b \in A_j$ can be realized as $\Vert b \Vert = \sup_{\phi \in S(A_j)}\vert \phi(b) \vert$. Interchanging suprema we then get
    \begin{align*}
        d_j(g,h) &= \sup \{ \sup_{\phi \in S(A_j)} \vert \phi(\alpha_g^{(j)} (a) - \alpha_h^{(j)}(a) ) \vert \mid a=a^*, K_j(a) \leq 1 \} \\ 
        &= \sup_{\phi \in S(A_j)} \mk_{K_j} (\phi \circ \alpha_g^{(j)}, \phi \circ \alpha_h^{(j)} ).
    \end{align*}
    Note that this metric is bi-invariant, so replacing $g$ and $h$ by their inverses leaves it unchanged.

    After passing to a subsequence, we therefore obtain, using \cite[Proposition 2.3 and Lemma 2.4]{Gelander2013}, a compact metric group $H$ along with almost inverse almost isometries between $H_j$ and $H$, with uniformly vanishing multiplicativity errors. Since the maps also have asymptotically dense images, the word estimates from \cite[Lemma 2.4]{ArzhantsevaCherix2024} imply that every fixed group identity satisfied by all the $H_j$ will also hold in $H$. In particular, $H$ is abelian and $h^p =e$ for all $h \in H$.

    For any $a \in \dom K_j$
    \begin{align*}
        \Vert \alpha_g^{(j)} (a) - \alpha_h^{(j)} (a) \Vert \leq d_j (g,h) K_j(a).
    \end{align*}
    So the actions have a common Lipschitz bound of one, and a common continuity modulus $t \mapsto t$. We may then apply \cite[Theorem 3.11]{LatremoliereCategories2019} to the subsequence and the limit $(A,K)$. Since every action is ergodic, we obtain a strongly continuous ergodic action of $H$ on $A$ by full quantum isometries. 
\end{proof}

\begin{lemma}\label{lemma:group-action-forces-AF}
    Let $p$ be a fixed prime, and suppose $H$ is a compact metrizable abelian group for which $h^p= e$ for every $h \in H$. If $H$ acts strongly continuously and ergodically on a nonzero unital $C^*$-algebra $A$, then
    \begin{align*}
        A \cong C((\bbZ_p)^r) \otimes B_{p,s}, \qquad r,s\in \bbN_0 \cup \{\infty\}
    \end{align*}
    where $B_{p,s}$ is as in Theorem \ref{thm:limits-are-AF}. In particular, $A$ must be AF. 
\end{lemma}

\begin{proof}
    Denote the action by $\alpha$,  let $G := H/\ker \alpha$, and $\Gamma = \widehat{G}$. As $G$ is compact, $\Gamma$ is discrete. The induced action by $G$ is faithful and ergodic. By \cite[§2.11]{OlesenPedersenTakesaki1980} each nonzero spectral subspace $A_x = \{a \mid \alpha_g(a) = x(g)a \text{ for all } g\in G\} $, $x \in \Gamma$, is generated by a unitary $u_x$, and their span is norm dense in $A$. 
    The nonzero spectral labels form a subgroup of $\Gamma$, since products and adjoints of these unitaries are nonzero. Its annihilator is the kernel of the action, so faithfulness and Pontryagin duality imply that every $x\in \Gamma$ occurs.
    Since $u_x u_y \in A_{x+y} = \bbC u_{x+y}$ there is a unique $\sigma(x,y) \in \bbT$ such that $u_xu_y = \sigma(x,y) u_{x+y}$. Associativity and the choice $u_0 = 1$ make $\sigma$ a normalized $2$-cocycle. 

    $G$ comes equipped with a normalized Haar measure, and we get $\int_G \alpha_g (a) dg = \tau(a)1$ for a faithful state $\tau$. Indeed, if $a \geq 0$ has zero average, then every continuous nonnegative function $g \mapsto \phi(\alpha_g(a))$ has zero integral, which implies $\phi(a)=0$ for every state and thus $a = 0$. In the associated faithful GNS representation $(\pi_\tau, H_\tau, \Omega_\tau)$ the vectors $u_x \Omega_\tau$ form an orthonormal basis, since $\tau(u_x^* u_y) = \delta_{x,y}$. We identify them with $\delta_x$ to obtain $\pi_\tau (u_x) \delta_y = \sigma(x,y) \delta_{x+y}$, $A \cong C^*_r(\Gamma, \sigma)\cong C^*(\Gamma, \sigma)$. The last isomorphism follows since $\Gamma$ is abelian and therefore amenable.  

    Because $G$ is compact metrizable, $\Gamma$ is countable. Moreover, $g^p = e$ implies $px = 0$ for every $x \in \Gamma$, and so it is an $\bbF_p$-vector space. The scalar commutators form a bicharacter with values in the pth roots of unity, so
    \begin{align*}
        u_x u_y u_x^* u_y^* = e^{2\pi i b(x,y)/p} \cdot 1
    \end{align*}
    defines an alternating bilinear form $b$. Let $R$ be its radical. We now rephase the unitaries indexed by a basis of $\Gamma$ to make their pth powers equal to $1$. 
    These generators present the full twisted group algebra as the universal algebra for their power and commutation relations. Thus \cite[Theorem 3.1]{ArvesonPrice2003} provides the isomorphism
    \begin{align*}
        A \cong Z \otimes B_{p,s}, \qquad Z = \overline{\operatorname{span}}\{ u_x \mid x \in R \}
    \end{align*}
    where $\dim_{\bbF_p} (\Gamma / R) = 2s$ in finite dimension and $s = \infty$ otherwise. The same argument covers finite $\Gamma$, with $B_{p,0} = \bbC$ in the rank zero case.

    To identify $Z$, choose a basis $(r_i)_i$ of $R$ and rephase $u_{r_i}$ to obtain commuting unitaries $z_i$ with $z_i^p = 1$. For any finite selection $r_1, \ldots , r_k$, the reduced monomials satisfy
    \begin{align*}
        z_1^{m_1} \cdots z_k^{m_k}= c_m u_{\sum_i m_i r_i}, \qquad 0\leq m_i < p \quad c_m \in \bbT.
    \end{align*}
    Their labels are distinct by linear independence of the $r_i$. Applying these operators to $\delta_0$ in the twisted regular representation gives nonzero scalar multiples of distinct basis vectors, so the monomials are linearly independent. Hence
    \begin{align*}
        C^*(z_1, \ldots, z_k) \cong C^*((\bbZ_p)^k). 
    \end{align*}
    Thus
    \begin{align*}
        Z \cong C^*(R) \cong C(\widehat{R}) \cong C((\bbZ_p)^r) \qquad r=\dim_{\bbF_p} R.
    \end{align*}
    This proves the asserted decomposition. $A$ being AF now follows from both factors in the decomposition being AF. 
\end{proof}

\begin{proof}[Proof of Theorem \ref{thm:limits-are-AF}]
    By Lemma \ref{lemma:limit-group-action} the limit algebra $A$ admits a strongly continuous and ergodic action of a compact metrizable abelian group $H$ with the prescribed property $h^p = e$ for all $h \in H$. By Lemma \ref{lemma:group-action-forces-AF} $A$ takes the form specified in the theorem statement. 
\end{proof}

\begin{lemma}\label{lemma:limit-length-spectral-triple}
    Let $p$ be a prime and let $G_i$ be either $\bbF_p$ with trivial cocycle or $\bbF_p^2$ with the standard nondegenerate Pauli cocycle. On $G = \oplus_{i\geq 1} G_i$ let $\sigma$ be the product cocycle, and use the length
    \begin{align*}
        \ell(g) = \sum_{i \geq 1} w_i 1_{g_i\neq e}, 
    \end{align*}
    where $w_i >0$ are such that $\sum_{i \geq 1} w_i^{-1} < \infty$. Let $A = C^*_r(G, \sigma)$ and $A_m = C^*_r(\oplus_{i=1}^m G_i, \sigma)$. Let $E_m \colon A \to A_m$ be the canonical conditional expectation which discards Fourier coefficients outside $\oplus_{i=1}^m G_i$. Then the spectral triple $(C^*_r(G, \sigma), \ell^2(G), D_\ell)$ is a metric spectral triple, its closed commutator seminorm restricts on each $A_m$ to the commutator seminorm for the restricted length, and, for $a \in \dom (K_\ell)$,
    \begin{align}\label{eq:tail-estimate}
        \Vert a - E_m (a) \Vert \leq (p^2 -1) \sum_{i >m} w_i^{-1} K_\ell(a).
    \end{align}
    As a result, if $G^{(m)} = \oplus_{i=1}^m G_i$, we have
    \begin{align*}
        \Lambda^{\rm dual} ((A, K_\ell), (A_m, K_{\ell\vert_{G^{(m)}}})) \leq (p^2 -1) \sum_{i > m} w_i^{-1}.
    \end{align*}
\end{lemma}

\begin{proof}
    Let $D_\ell$ be the self-adjoint multiplication operator on $\ell^2(G)$ with
    \begin{align*}
        \dom D_\ell = \{ \xi \in \ell^2(G) \mid \sum_{g\in G} \ell(g)^2 \vert \xi(g) \vert^2 < \infty \}.
    \end{align*}
    Use the maximal commutator domain and lower semicontinuous seminorm as in Section \ref{sec:preliminaries}.  Note that as $w_i \to \infty$, the length balls are finite and $D_\ell$ has compact resolvent. Moreover, the $*$-algebra $C_c(G, \sigma)$  is dense in $C^*(G, \sigma)$ and consists of elements for which the commutators are bounded since $\Vert [D_\ell, \lambda_g^\sigma ]\Vert \leq \ell(g)$. 

    Now fix an integer $k \geq 0$, and set
    \begin{align*}
        G^{(k)} := \bigoplus_{i \leq k} G_i, \qquad G^{(>k)} := \bigoplus_{i > k} G_i. 
    \end{align*}
    Let $A_{>k}$ and $D_{>k}$ denote the twisted subgroup algebra and multiplication operator for the restricted length, and define $D_k$ similarly. The product cocycle and additive length give
    \begin{align*}
        \ell^2(G) &= \ell^2(G^{(k)}) \otimes \ell^2(G^{(>k)}), \qquad A=A_k \otimes  A_{>k} \\
        D_\ell &= D_k \otimes I + I \otimes D_{>k}
    \end{align*}
    The operator equality includes domains, since $D_k$ is bounded and both sides have domain $\dom (I \otimes D_{>k})$. Thus
    \begin{align*}
        [D_\ell, a \otimes I] = [D_k, a] \otimes I, \qquad K_\ell(a) = K_{\ell\vert_{G^{(k)}}}(a), \qquad a\in A_k.
    \end{align*}
    Moreover, $E_k = \id \otimes \tau_{>k}$, where $\tau_{>k}$ is the canonical trace, is a completely positive retraction. 
    For $V_k \xi = \xi \otimes \delta_e$ we have $E_k(a) = V_k^* a V_k$ and $D_\ell V_k = V_k D_k$, so, for $a \in \dom(K_\ell)$
    \begin{align*}
        [D_k,E_k(a)]=V_k^*[D_\ell,a]V_k,\qquad K_\ell(E_k(a))\leq K_\ell(a).
    \end{align*}
    Also $E_k(a) \to a$ in norm, first for finite Fourier polynomials and then for every $a \in A$ by contractivity.
    Let $Q_i$ be the conditional expectation which applies the canonical trace to the $i$th factor. These maps commute, and averaging the dual action on the factor shows that each is Lip-contractive. For a finite Fourier polynomial $a = \sum_{u \in G_i} a_u \otimes \lambda_u^{\sigma_i}$, write $D_{\rm rest}$ for the multiplication operator on the remaining factors. The two-corner argument in the proof of Proposition \ref{prop:concatentation-qudit-collapse} gives
    \begin{align*}
        B_{\pm} = [D_{\rm rest} , a_u] \pm w_i a_u, \qquad \Vert B_{\pm } \Vert \leq K_\ell(a), \quad u\neq e.
    \end{align*}
    Subtracting yields $\Vert a_u \Vert \leq w_i^{-1}K_\ell(a)$, and therefore 
    \begin{align*}
        \|a-Q_i(a)\|\leq\sum_{u\ne e}\|a_u\|\leq\frac{|G_i|-1}{w_i}K_\ell(a).
    \end{align*}
    For $a \in A_N$ and $N \geq k$ we have $E_k(a) = Q_{k+1} \cdots Q_N(a)$. Applying the preceding estimate successively, using Lip-contractivity at each step, gives 
    \begin{align*}
        \|a-E_k(a)\| \leq \sum_{i=k+1}^N\frac{|G_i|-1}{w_i}K_\ell(a)  \le(p^2-1)\left(\sum_{i=k+1}^Nw_i^{-1} \right) K_\ell(a).
    \end{align*}
    For arbitrary $a \in \dom (K_\ell)$, apply this inequality to $E_N(a) \in A_N$. Since 
    \begin{align*}
        K_\ell(E_N(a))\le K_\ell(a),\qquad E_k E_N(a)=E_k(a)
    \end{align*}
    for $N \geq k$, letting $N \to \infty$ and using $E_N(a) \to a$ in norm proves \eqref{eq:tail-estimate}.

    For $k=0$, the algebra $A_0$ is $\bbC$, and $E_0(a) = \tau(a) 1$. As such, the tail estimate provides a scalar kernel, and a bound on the self-adjoint Lip-unit ball centered in $\tau$. Its image under each $E_k$ is a bounded subset of a finite-dimensional algebra $A_k$, and is therefore totally bounded. Then \eqref{eq:tail-estimate} approximates the full Lip-unit ball uniformly as $k\to \infty$, proving that the $\tau$-centered Lip-unit ball is totally bounded and so the spectral triple is metric. Lemma \ref{lemma:cp-approximation-map} now provides the propinquity bound. 
    Note that for a finite list the same argument would just use finite sums, and once the tail vanishes the pairs will coincide. 
\end{proof}

Recall the codes $\calF^{(j)}$ from Lemma \ref{lemma:513-protector-code} above. Lemma \ref{lemma:limit-length-spectral-triple} then allows us to use them to realize specific UHF algebras as limit points in $\overline{\calQ_p}$.
\begin{theorem}\label{theorem:CAR-limit}
    For $\calC_m = \otimes_{j=1}^m \calF^{(j)}$ the code parameters and length are given by
    \begin{align*}
        [[n_m, k_m, d_m ]]_p &= \left[\left[ \frac{5^{m+1}-5}{4}  , m , 3\right]\right]_p \\
        \ell_{\calC_m} (g) &= \sum_{j=1}^m 3^j 1_{g_j \neq e}
    \end{align*}
    for $g \in \oplus_{j=1}^m \bbF_p^2$. Let $G_\infty = \oplus_{j\geq 1} \bbF_p^2$ with its product cocycle, and set
    \begin{align*}
        A_\infty := C^*_r (G_\infty , \sigma) \cong \bigotimes_{j\geq 1} M_p (\bbC), \qquad D_\infty \delta_g := \left( \sum_{j\geq 1} 3^j 1_{g_j \neq e} \right) \delta_g.
    \end{align*}
    This yields a metric spectral triple, and its associated pair satisfies
    \begin{align*}
        \Lambda^{\rm dual} ((A_{\calC_m} , K_{\calC_m}) , (A_{\infty} , K_{D_\infty})) \leq \frac{p^2 -1}{2} \cdot 3^{-m}.
    \end{align*}
    In particular, for $p=2$ the limit is a metric spectral triple on the CAR algebra.
\end{theorem}

\begin{proof}
    The tensor factors act on disjoint physical sites, so their stabilizer and logical spaces are direct sums, and the quotient Hamming lengths add as in Proposition \ref{prop:concatentation-qudit-collapse}. 
    Lemma \ref{lemma:513-protector-code} gives $\dim \calF^{(j)} = p$ and length $3^j$ for every nonzero logical label in the $j$th factor. It follows that $k_m = m$, $d_m = 3$, and $n_m = \sum_{j=1}^m 5^j = (5^{m+1} -5)/4$. The resulting spectral triples are exactly the first $m$ tensor factors of the spectral triple on $A_\infty$. Now apply Lemma \ref{lemma:limit-length-spectral-triple} with $w_j = 3^j$, and use $\sum_{j>m} 3^{-j} = 3^{-m}/2$ to finish the proof. 
\end{proof}

The preceding constructions can be combined to realize all the algebras in Theorem \ref{thm:limits-are-AF}, using ordinary Hamming weights throughout.
\begin{theorem}\label{thm:AF-boundary-realization}
 For every $r,s\in\{0,1,2,\ldots,\infty\}$, the algebra
 $C((\bbZ_p)^r)\otimes B_{p,s}$ admits a metric spectral triple whose associated pair belongs to $\overline{\mathcal{Q}_p}$.
\end{theorem}
\begin{proof}
   Assume first that $(r,s) \neq(0,0)$.  Choose a finite initial segment $I\subseteq\bbN$, or $I=\bbN$, partitioned as
   \begin{align*}
        I=I_{\mathrm{cl}}\sqcup I_{\mathrm{q}},  \qquad |I_{\mathrm{cl}}|=r,\qquad |I_{\mathrm{q}}|=s.
   \end{align*}
   Choose coordinates on $\calL_{\calF^{(i)}}\cong\bbF_p^2$, and logical generators with associated cocycle
   \begin{align*}
       \sigma_W ((a,b) , (a', b')) = \zeta^{ba'}. 
   \end{align*}
   Put $J=\{0\}\times\bbF_p$ and define
   \begin{align*}
       G_i=\begin{cases}
            J,&i\in I_{\mathrm{cl}},\\
            \bbF_p^2,&i\in I_{\mathrm{q}},
            \end{cases}
 \qquad
 A_i=\Cred(G_i,\sigma_{\mathrm W}|_{G_i\times G_i})\cong
 \begin{cases}
 C(\bbZ_p),&i\in I_{\mathrm{cl}},\\
 M_p(\bbC),&i\in I_{\mathrm{q}}.
 \end{cases}
   \end{align*}     
Indeed, the cocycle is trivial on $J$, whose labels $(0,b)$ act as powers of the logical generator $\overline Z_i$. Thus the commutative  choice is the subalgebra $C^*(\overline Z_i)$ of the full logical algebra.

Now set
\begin{align*}
    G=\bigoplus_{i\in I}G_i,\qquad \sigma=\prod_{i\in I}\sigma_{\mathrm W}|_{G_i\times G_i}, \qquad \ell(g)=\sum_{i\in I}w_i\mathbf1_{\{g_i\ne0\}}, \quad w_i=3^i.
\end{align*}
   Then $A=\Cred(G,\sigma)\cong\bigotimes_{i\in I}A_i \cong C((\bbZ_p)^r)\otimes B_{p,s}.$
   Lemma~\ref{lemma:limit-length-spectral-triple} supplies its metric spectral triple and the finite factor algebras $B_h=\Cred(\bigoplus_{i\in I,\ i\le h}G_i,\sigma)$, with their regular seminorms, satisfying
   \begin{align}\label{eq:boundary-tail}
       \Lambda^{\rm dual}((A,K_A),(B_h,K_{B_h})) \le(p^2-1)\sum_{\substack{i\in I\\i>h}}w_i^{-1}.
   \end{align}
   It thus remains to approximate $B_h$ by the logical algebra of a code. 

   Fix $h \in I$ and $N \geq 1$. Choose  unitaries  $V_i \colon \bbC^p \to \calF^{(i)}$ implementing the chosen $*$-isomorphisms $M_p(\bbC) \cong B(\calF^{(i)})$, and define
   \begin{align*}
       \mathcal{D}_{i,N} = \begin{cases}
           V_i^{\otimes N} (\mathcal{R}_N), & i\in I_{\rm cl} \\
           \calF^{(i)}, & i\in I_{\rm q},
       \end{cases}
       \qquad \calC_{h,N} = \bigotimes_{i=1}^h \mathcal{D}_{i,N}
   \end{align*}
   where $\mathcal{R}_N$ is the repetition code from Example \ref{example:repetition-code}. Each factor has dimension $p$. 

   Since every nonzero logical class of $\calF^{(i)}$ has length $w_i$, the same argument as in Lemma \ref{lemma:513-protector-code} multiplies the quotient lengths of $\mathcal{R}_N$ by $w_i$. We obtain
   \begin{align*}
       \ell_{i,N} (a,b) &=w_i\mathbf1_{\{a=0,\ b\ne0\}}
		+Nw_i\mathbf1_{\{a\ne0\}},
		\qquad &&i\in I_{\mathrm{cl}}, \\
        \ell_{i,N}(a,b)&=w_i\mathbf1_{\{(a,b)\ne(0,0)\}} &&i \in I_{\rm q}.
   \end{align*}
    As in Proposition~\ref{prop:concatentation-qudit-collapse}, these lengths add under tensor products.
   Thus $A_{\calC_{h,N}}\cong M_{p^h}(\bbC)$, and $B_h$ is its subgroup algebra corresponding to $\bigoplus_{i=1}^hG_i$.

   	We next verify exact seminorm restriction. On the regular Hilbert space of $B_h$, write
    \begin{align*}
        D_h=\sum_{i=1}^h T_i, \qquad T_i=w_i( I-P_{\bbC \delta_0}),
    \end{align*}
    where $T_i$ acts on the $i$th tensor factor.
    The ambient regular representation decomposes over cosets of $\bigoplus_{i=1}^hG_i$. By the calculation in Example \ref{example:repetition-code}, a diagonal unitary identifies the subgroup representation on each coset with the regular representation of $B_h$.
    Under this identification, the ambient Dirac operator is $D_h-\sum_{i\in O}T_i+cI,$  where $O$ consists of the indices $i \in I_{\rm cl}$ for which the coset has nonzero  first coordinate, and $c$ is a scalar.

    For $i \in I_{\rm cl}$, let $U_i$ be the cyclic shift on $\ell^2(G_i)$, acting identically on the other factors. It commutes with $B_h$, and
    \begin{align*}
        \frac1p\sum_{b=0}^{p-1}U_i^bT_iU_i^{-b} =w_i\left(1-\frac1p\right)I.
    \end{align*}
    Consequently, for $a=a^*\in B_h$,
    \begin{align*}
        [D_h-T_i,a] =\frac1p\sum_{b=0}^{p-1}U_i^b[D_h,a]U_i^{-b},\qquad\|[D_h-T_i,a]\|\le\|[D_h,a]\|.
    \end{align*}
    Iterating bounds every coset commutator by $\|[D_h,a]\|$. The identity coset has $O=\varnothing$ and gives exactly this commutator. Therefore $K_{\calC_{h,N}}|_{(B_h)_{\rm sa}}=K_{B_h}$. 

    The expectation onto $B_h$ is Lip-contractive by Lemma~\ref{lemma:fourier-truncation}.
    Every discarded label has nonzero first coordinate for some $i \in I_{\rm cl}$, and hence length at least $N\min_{i\le h}w_i$. There are at most $p^{2h}$ labels.
    The Fourier estimate and
	Lemma \ref{lemma:cp-approximation-map} therefore give
    \begin{align}\label{eq:boundary-collapse}
        \Lambda^{\rm dual} \bigl((A_{\calC_{h,N}},K_{\calC_{h,N}}),(B_h,K_{B_h})\bigr) \le\frac{p^h}{N\min_{i\le h}w_i}.
    \end{align}

    	To obtain a single convergent sequence, choose
        \begin{align*}
            h_m=\begin{cases}
		          \min\{m,|I|\},&I\text{ finite},\\
		          m,&I=\bbN,
	               \end{cases}
	           \qquad N_m=mp^{h_m}.
        \end{align*}
        Since $\min_{i\le h_m}w_i=3$, the triangle inequality and \eqref{eq:boundary-tail} and \eqref{eq:boundary-collapse} yield
        \begin{align*}
            \Lambda^{\rm dual} \bigl((A_{\calC_{h_m,N_m}},K_{\calC_{h_m,N_m}}),(A,K_A)\bigr) \le \frac{1}{3m} +(p^2-1)\sum_{\substack{i\in I\\i>h_m}}3^{-i} \longrightarrow 0.
        \end{align*}
        All source lengths arise from ordinary physical Hamming weights.

        Finally, if $r=s=0$, then $A=\bbC$, and Proposition \ref{prop:concatentation-qudit-collapse} supplies the required approximation.
\end{proof}

\begin{corollary}\label{cor:AF-list-realized}
    The underlying $C^*$-algebras occurring in $\overline{\calQ_p}$ are exactly the algebras in \eqref{eq:boundary-algebra-classification}. In particular, every algebra underlying a metric boundary point is AF. Each listed algebra which is not isomorphic to $M_{p^k}(\bbC)$ for some $k \geq 1$ occurs as the algebra of an added boundary point. 
\end{corollary}

\begin{proof}
    Any algebra in the limit must be of the listed type by Theorem \ref{thm:limits-are-AF}. We know every algebra in the list is realized by Theorem \ref{thm:AF-boundary-realization}. Lastly, a pair whose algebra is not a full matrix algebra of the indicated size cannot belong to the finite-code image, because zero propinquity requires a $*$-isomorphism.
\end{proof}

\subsection{Error correction as a local property}

Returning to codes and thus finite-dimensional $C^*$-algebras of $\calQ_p$, we show that error correction of certain collections of possible errors is a local property which is built into the topology. Note however, that by virtue of how we have defined the distance between codes, correction of arbitrary prescribed families of errors cannot be a property of a point in $\calQ_p$.  Indeed, by simple interchanging of sites the codes $|0\rangle \otimes \bbC^p$ and $\bbC^p \otimes |0 \rangle$ have distance zero. However, the error space $\operatorname{span}\{ I, Z_{p,1} \}$ is correctable on the first code since it acts by scalars, but fails the Knill--Laflamme criterion on the second. 
Under the hypotheses of Theorem \ref{thm:local-error-correction}, correction of all errors up to a prescribed weight is a property of a point and is locally constant. 

The argument will make use of an intrinsic quantity defined on the self-adjoint part of the algebra. If $(A,K)$ is a Leibniz quantum compact metric space with tracial state $\tau$, define
\begin{align}\label{eq:def-kappa-intrinsic}
    \kappa (A,K, \tau ) = \inf \{ K(a) \mid a=a^*, \tau(a) = 0, \tau(a^2) = 1 \}. 
\end{align}
If $A$ is finite-dimensional, $A \not \cong \bbC$, and $\tau$ is a faithful tracial state, this infimum is finite and attained. Note that the quantity is a priori larger than the similar quantity defined in Definition \ref{def:code distance}. However, for Pauli stabilizer codes and $p \in \{2,3,5\}$ they will agree, see Lemma \ref{lemma:kappa-versus-distance} below. 

\begin{proposition}\label{prop:kappa-continuous}
    Let $(A_j, K_j)$ be Leibniz quantum compact metric spaces for which $A_j$ is a nontrivial full matrix algebra and $\tau_j$ the normalized trace. Suppose $(A_j, K_j)$ converges to $(A,K)$ in dual propinquity, where $A$ is itself a nontrivial full matrix algebra, and where $\tau$ is the normalized trace. Then
    \begin{align*}
        \kappa(A_j, K_j, \tau_j) \longrightarrow \kappa(A, K, \tau).
    \end{align*}
\end{proposition}

\begin{proof}
    By definition of dual propinquity, we choose tunnels $(C_j, M_j, \rho_j, \pi_j)$ between $A_j$ and $A$ with extents $\varepsilon_j \to 0$. Again by definition we may pick $\varphi_j \in S(A)$ with $\mk_{M_j} (\tau_j \circ \rho_j , \varphi_j \circ \pi_j) \leq \varepsilon_j$. For fixed $a,b \in A_{\rm sa}$ pick self-adjoint lifts $x_j, y_j$ through $\pi_j$ with $M_j(x_j) \leq K(a) + 1/j$ and $M_j(y_j) \leq K(b) + 1/j$. Now \eqref{eq:lift-norm-bounds} shows that the lifts are uniformly bounded in norm, and  so \eqref{eq:Leibniz-condition} bounds $M_j(\{x_j,y_j\})$ uniformly.
    Since $\tau_j (\{ \rho_j(x_j) , \rho_j(y_j)\}) =0$ we obtain $|\varphi_j(\{a,b\})|\le\varepsilon_j M_j(\{x_j,y_j\})\to 0$. Every cluster state of $\varphi_j$ is therefore tracial. 
    The trace on $A$ is unique, so $\varphi_j\to\tau$ weak$^*$ and in $\mk_K$. Consequently, for every  $x \in \dom (M_j)$,
    \begin{align}\label{eq:rj-to-zero}
        |\tau_j(\rho_j(x))-\tau(\pi_j(x))| \le r_j M_j(x),  \quad \text{where} \quad r_j=\varepsilon_j + \mk_K(\varphi_j,\tau) \longrightarrow 0.
    \end{align}
    Let $a$ attain the minimum in \eqref{eq:def-kappa-intrinsic}, and choose a self-adjoint lift $x_j$ satisfying $\pi_j(x_j)=a$ and $M_j(x_j)\le K(a)+1/j$. The norms of $x_j$ are bounded, and $M_j(x_j^2)\le 2 \|x_j\| M_j(x_j)$. Apply \eqref{eq:rj-to-zero} to $x_j$ and $x_j^2$. Then with $b_j = \rho_j(x_j)$ we have $\tau_j(b_j) \to 0$ and $\tau_j(b_j^2) \to 1$. Thus, for all sufficiently large $j$ we may define
    \begin{align*}
        c_j= \frac{b_j-\tau_j(b_j)1_{A_j}} {\sqrt{\tau_j(b_j^2)-\tau_j(b_j)^2}}.
    \end{align*}
    Then $\tau_j(c_j)=0$ and $\tau_j(c_j^2)=1$. Since subtracting a scalar does not change $K_j$, we obtain
    \begin{align*}
        \kappa(A_j,K_j,\tau_j) \le K_j(c_j)\le \frac{K(a)+1/j} {\sqrt{\tau_j(b_j^2)-\tau_j(b_j)^2}} \longrightarrow K(a) =\kappa(A,K,\tau).
    \end{align*}

    For the reverse inequality, choose self-adjoint minimizers $a_j\in A_j$, so that $ \tau_j(a_j)=0$, and $\tau_j(a_j^2)=1$, $K_j(a_j)=\kappa(A_j,K_j,\tau_j)$.  
    The preceding upper bound shows that $K_j(a_j)$ is bounded. Since $\tau_j(a_j)=0$, we get by definition of extent
    \begin{align*}
        \|a_j\| = \sup_{\phi \in S(A_j)} \vert \phi(a_j) - \tau_j(a_j)\vert \le \operatorname{diam}(S(A_j),\mk_{K_j})K_j(a_j) \le\bigl(\operatorname{diam}(S(A),\mk_K)+2\varepsilon_j\bigr) K_j(a_j).
    \end{align*}
    Thus the norms $\|a_j\|$ are bounded as well.

    Choose self-adjoint lifts $y_j\in C_j$ satisfying $\rho_j(y_j)=a_j$ and  $M_j(y_j)\le K_j(a_j)+1/j$. Then $\|y_j\|\le\|a_j\|+\varepsilon_jM_j(y_j)$,  so $y_j$ is bounded. Moreover
    \begin{align*}
        K(\pi_j(y_j))\le M_j(y_j),\qquad M_j(y_j^2)\le2\|y_j\|M_j(y_j).
    \end{align*}
    Applying \eqref{eq:rj-to-zero} to $y_j$ and $y_j^2$, and using that the endpoint maps
are $*$-homomorphisms, yields
\begin{align*}
     |\tau(\pi_j(y_j))| &\le r_j M_j(y_j)\longrightarrow 0,\\ 
     |\tau(\pi_j(y_j)^2)-1|  &\le 2r_j\|y_j\|M_j(y_j)\longrightarrow0.
\end{align*}
Put $b_j = \pi_j(y_j)$. For sufficiently large $j$, centering and normalizing $b_j$ gives
\begin{align*}
    \kappa(A, K ,\tau) \leq K \left( \frac{b_j - \tau(b_j) 1_A}{\sqrt{\tau(b_j^2) - \tau(b_j)^2}}\right) \leq \frac{K_j(a_j) + 1/j}{\sqrt{\tau(b_j^2) - \tau(b_j)^2}}. 
\end{align*}
The denominator tends to one, so $\kappa(A,K,\tau) \leq \liminf_j \kappa(A_j, K_j, \tau_j)$. 
Together with the previously established upper limit inequality, this establishes the asserted convergence. 
\end{proof}

We next compare $\kappa(A_\calC,K_\calC,\tau)$ with the minimum code distance. Recall that $\calL_\calC$ is an
$\bbF_p$-vector space. For $t\in\bbF_p\setminus\{0\}$, multiplication by $t$ preserves the support of every physical Pauli label and bijects $S(\calC)$ onto itself. It therefore maps the representatives of $v$ bijectively onto those of $tv$, preserving their Hamming weights. Taking minima as in
\eqref{eqn:length function on logical algebra} gives
\begin{align}\label{eq:length-invariance-Fp}
 \ell_\calC(tv)=\ell_\calC(v)  
\end{align}
for $v\in\calL_\calC$, $t\in\bbF_p\setminus\{0\}$.

\begin{lemma}\label{lemma:kappa-versus-distance}
    Let $\calC$ be a prime-$p$ Hamming stabilizer code with $\dim \calC > 1$ and distance $d_\calC = \Delta_\calC$. Then 
    \begin{align*}
        \kappa(A_{\calC},K_{\calC},\tau)=d_{\calC} \qquad\text{if }p\in\{2,3,5\}.
    \end{align*}
     For every odd prime one also has
     \begin{align}\label{eq:general-prime-bound}
         d_{\calC}\le\kappa(A_{\calC},K_{\calC},\tau)\le c_p d_{\calC}, \qquad c_p=\frac{p}{\sqrt{p^2-1}}.
     \end{align}
\end{lemma}

\begin{proof}
    Write $d = d_\calC$. 
    For $a=a^*=\sum_x a_x\lambda_x^{\sigma_\gamma}$, orthogonality $\tau((\lambda_x^{\sigma_\gamma})^*\lambda_y^{\sigma_\gamma}) =\delta_{x,y}$ and $\tau(a)=a_0$ give $\tau(a^2)-\tau(a)^2=\sum_{x\ne0}|a_x|^2$ by Parseval's identity. Combining this with \eqref{eq:lower-2-bound} gives
    \begin{align*}
        K_\calC (a)^2 \geq  \sum_x |a_x|^2\ell_\calC(x)^2
 \geq d^2\bigl(\tau(a^2)-\tau(a)^2\bigr)
    \end{align*}
    for $a=a^*$, thus proving the common lower bound. 
    Now choose a shortest nonzero logical label $v$ and choose a scalar multiple $U$ of $\lambda_v^{\sigma_\gamma}$ such that $U^p = I$. 
    On a coset $x+\langle v\rangle$, $U$ cyclically permutes the vectors $\delta_{x+rv}$, $r\in\bbF_p$, up to unimodular coefficients. Since $U^p=I$, their product is one, so a diagonal unitary change of basis removes these coefficients while leaving the diagonal Dirac operator unchanged. Its eigenvalues on this subspace satisfy $|\ell_\calC(x+rv)-\ell_\calC(x+sv)| \le \ell_\calC((r-s)v)\le d $, for $r,s\in\bbF_p$ by the triangle inequality and \eqref{eq:length-invariance-Fp}. 

    We now cover the various cases for the value $p$. 
    For $p=2$, the unitary $U$ is self-adjoint and has trace zero and trace variance one. Since $K_\calC(U)=d$, it attains the lower bound.

    For $p=3$, put $P=(I+U+U^*)/3$ and $a=(3P-I)/\sqrt2$. Then $P$ is a projection with $\tau(P)=1/3$, so $\tau(a)=0$ and $\tau(a^2)=1$. On each coset block, $P$ projects onto $\xi = (1,1,1)/\sqrt{3}$.
    Writing $h=(I-P)D\xi$, the decomposition $\bbC\xi\oplus\xi^\perp$ gives
    \begin{align*}
        [D,P]=\begin{pmatrix}0&-h^*\\h&0\end{pmatrix}, \qquad \|[D,P]\|=\|h\|.
    \end{align*}
    Thus, for the Dirac values $r_0,r_1,r_2$ on this block,
   \begin{align*}
    \|[D,P]\|^2 &=\|D\xi\|^2-|\langle D\xi,\xi\rangle|^2 =\frac13\sum_{i=0}^2r_i^2 -\left(\frac13\sum_{i=0}^2r_i\right)^2 \le\frac{2d^2}{9}.
\end{align*}
Indeed, translating the $r_i$ into $[0,d]$ preserves their variance. This convex quadratic attains its maximum at a vertex, where its value is either zero or $2d^2/9$. Since $[D,a]=(3/\sqrt2)[D,P]$, taking the maximum over coset blocks yields $K_\calC(a)\le\frac{3}{\sqrt2}\sqrt{\frac{2d^2}{9}}=d$. The common lower bound therefore gives
$\kappa(A_\calC,K_\calC,\tau)=d$.

For $p=5$, put
\begin{align*}
    c_1=\frac{e^{i\pi/5}}2,\qquad c_2=\frac{e^{2i\pi/5}}2,\qquad a=c_1U+\overline{c_1}U^* +c_2U^2+\overline{c_2}U^{*2}.
\end{align*}
Using $U^5=I$, multiplication gives
\begin{align*}
    a^2=I+\frac12\left( c_1U+\overline{c_1}U^*-c_2U^2-\overline{c_2}(U^*)^2\right).
\end{align*}
Since $\tau(U^j)=0$ for $1\le j\le4$, we have $\tau(a)=0$ and $\tau(a^2)=1$.

On each cyclic coset block, $a$ has zero diagonal and off-diagonal entries of modulus $1/2$, while $a^2$ has diagonal entries one and off-diagonal entries of modulus $1/4$. For a coordinate projection $Q$, the off-diagonal blocks of $[Q,a]$ are $T=Qa(I-Q)$ and $-T^*$, so
\begin{align*}
    \|[Q,a]\|^2=\|TT^*\| =\|Qa^2Q-(QaQ)^2\|.
\end{align*}
For rank-one $Q$, the last operator is $Q$. For rank-two $Q$, $(QaQ)^2=Q/4$, so on $\operatorname{ran}Q$ the last operator has matrix
\begin{align*}
    \begin{pmatrix} 3/4 & z\\ \overline z & 3/4
 \end{pmatrix},
 \qquad |z|=1/4.
\end{align*}
Its eigenvalues are $1$ and $1/2$, hence its norm is one.
Using $[I-Q,a]=-[Q,a]$, we conclude that $\|[Q,a]\|\le1$ for every coordinate projection $Q$.

Order the Dirac values as $r_1\le\cdots\le r_5$, and let $Q_j$ project onto the coordinates corresponding to $r_{j+1},\ldots,r_5$. Then
\begin{align*}
     D=r_1 I+\sum_{j=1}^4(r_{j+1}-r_j) Q_j,
\end{align*}
and consequently
\begin{align*}
     \|[D,a]\| \le\sum_{j=1}^4(r_{j+1}-r_j)\|[Q_j,a]\| \le r_5-r_1\le d.
\end{align*}
Taking the maximum over coset blocks gives $K_\calC(a)\le d$. The common lower bound therefore yields $\kappa(A_\calC,K_\calC,\tau)=d$.

 Finally, let $p$ be any odd prime and put $m=(p-1)/2$.
 With $\zeta = e^{2\pi i/p}$, spectral projections of $U$ are
 \begin{align*}
     P_r=\frac{1}{p}\sum_{j=0}^{p-1}\zeta^{-rj}U^j,
 \qquad 0\le r<p,
 \end{align*}
 by functional calculus. Since $\tau(U^j)=0$ for $1\le j<p$, each satisfies $\tau(P_r)=1/p$. Thus
 \begin{align*}
     P=\sum_{r=0}^{m-1}P_r,\qquad  a=\frac{P-(m/p)I}{\sqrt{(m/p)(1-m/p)}}
 \end{align*}
 satisfy $P^2=P=P^*$, $\tau(a)=0$, and $\tau(a^2)=1$.

 On each coset block, let $c$ be the midpoint of the smallest and largest Dirac values. Their difference is at most $d$, so $\|D-cI\|\le d/2$. Relative to $\operatorname{ran}P\oplus\ker P$, we have
 \begin{align*}
     [D,P]=
 \begin{pmatrix}0&-T\\T^*&0\end{pmatrix},
 \qquad T=P(D-cI)(I-P).
 \end{align*}
 As such, $\|[D,P]\|=\|T\|\le\|D-cI\|\le d/2$. Taking the maximum over coset blocks gives $K_\calC(P)\le d/2$. Therefore
 \begin{align*}
     K_\calC(a) =\frac{K_\calC(P)}{\sqrt{(m/p)(1-m/p)}} \le\frac{d/2}{\sqrt{(p^2-1)/(4p^2)}} =\frac{p}{\sqrt{p^2-1}}\cdot d.
 \end{align*}
 Since $a$ is admissible in the definition of $\kappa$, the common lower bound yields
 \begin{align*}
     d\le\kappa(A_\calC,K_\calC,\tau) \le\frac{p}{\sqrt{p^2-1}}\cdot d.
 \end{align*}
\end{proof}

\begin{theorem}
\label{thm:local-error-correction}
 Let $\calC$ be a prime-$p$ Hamming stabilizer code with $\dim \calC > 1$ and minimum distance $d$. Assume either
 \begin{enumerate}
  \item $p\in\{2,3,5\}$, or
  \item $p\ge7$ and $d(c_p-1)<1$, where $c_p=p/\sqrt{p^2-1}$.
 \end{enumerate}
 Then there is $\varepsilon>0$ such that every prime-$p$ Hamming stabilizer code $\calC'$ with $\dim \calC' >1$ and $d_p(\calC,\calC')<\varepsilon$ has minimum distance $d$. In particular, for each integer $t\ge0$, $\calC'$ corrects the entire error space of weight at most $t$ if and only if $\calC$ does. Neither physical nor logical dimensions are fixed in this neighborhood.

 Condition (2) holds in particular whenever $d\le2p^2-2$. For $p\in\{2,3,5\}$ and fixed $t$, the set of finite-code points correcting all errors of weight at most $t$ is both open and closed in $\mathcal{Q}_p$.
\end{theorem}

\begin{proof}
    Suppose $d_p (\calC_j, \calC) \to 0$. Write $d_j$ for the distance associated with the $\calC_j$, and write $\kappa_j, \kappa$ for the quantities in \eqref{eq:def-kappa-intrinsic}. By Proposition \ref{prop:kappa-continuous} we have $\kappa_j \to \kappa$. For $p \in \{2,3,5\}$, Lemma \ref{lemma:kappa-versus-distance} gives $d_j \to d$. For any prime larger than this it gives
    \begin{align*}
        \frac{d}{c_p} \leq \frac{\kappa}{c_p} \leq \liminf_j d_j \leq \limsup_j d_j \leq \kappa \leq c_p d.
    \end{align*}
    If $d(c_p -1) <1$ the outer endpoints lie in the interval $(d-1,d+1)$, and the fact that $d_j$ and $d$ are integers forces $d_j=d$ eventually. The asserted statement about $\varepsilon$-neighborhoods follows. The sufficient integer bound follows from
 \begin{align*}
  \frac1{c_p-1}  =p\sqrt{p^2-1}+p^2-1>2p^2-2.
 \end{align*}

 It remains to connect distance to correction. Let $\calE_t=\Span\{\pi_{p,n}(x):\ell_{p,n}(x)\le t\}$. Products of errors satisfy $\Span(\calE_t^*\calE_t)=\calE_{\min\{2t,n\}}=\calE_{2t}$, where $\calE_s=\calE_n$ for $s\ge n$, since every Pauli
operator in $B((\bbC^p)^{\otimes n})$ has Hamming weight at most $n$. 
For the reverse inclusion, partition the support of a Pauli operator of weight at most $2t$ into two subsets of size at most $t$. The Knill--Laflamme criterion says that $\calE_t$ is correctable exactly when $P_\calC \calE_{2t}P_\calC \subseteq \bbC P_\calC$. A Pauli label outside $S(\calC)^\perp$ compresses to zero, a stabilizer label compresses to a scalar, and every label in $S(\calC)^\perp\setminus S(\calC)$ gives a nonscalar operator on $\calC$. Thus $\calE_t$ is correctable on $\calC$ if and only if $d_\calC > 2t$. No lower bound on the weight of a nonidentity stabilizer element is required. Local constancy of $d$ proves the correction statement, and in the case $p\in \{2,3,5\}$, the claim about open and closed sets. 
\end{proof}

\vspace{1cm} \noindent {\bf Use of AI assistance:} The writing of this article was assisted by ChatGPT-6.0 Astra, a large language model made by OpenAI. All AI-produced suggestions have been verified, and where necessary corrected, by the authors. 

\printbibliography

@misc{Eidesen2026,
  title = {Projective error models: {S}tabilizer codes, {C}lifford codes, and weak stabilizer codes},
  author = {Jonas Eidesen},
  year = {2026},
  eprint = {2506.01843v3},
  archivePrefix = {arxiv}
}

@article{KapshikarKundu2023,
  author={Kapshikar, Upendra and Kundu, Srijita},
  journal={IEEE Transactions on Information Theory}, 
  title={On the Hardness of the Minimum Distance Problem of Quantum Codes}, 
  year={2023},
  volume={69},
  number={10},
  pages={6293-6302},
}

@article {ChienWaldron17,
    AUTHOR = {Chien, Tuan-Yow and Waldron, Shayne},
     TITLE = {Nice error frames, canonical abstract error groups and the
              construction of {SIC}s},
   JOURNAL = {Linear Algebra Appl.},
  FJOURNAL = {Linear Algebra and its Applications},
    VOLUME = {516},
      YEAR = {2017},
     PAGES = {93--117},
}

@article {KlappeneckerRottelerI02,
    AUTHOR = {Klappenecker, Andreas and R\"otteler, Martin},
     TITLE = {Beyond stabilizer codes. {I}. {N}ice error bases},
   JOURNAL = {IEEE Trans. Inform. Theory},
  FJOURNAL = {Institute of Electrical and Electronics Engineers.
              Transactions on Information Theory},
    VOLUME = {48},
      YEAR = {2002},
    NUMBER = {8},
     PAGES = {2392--2395},
}

@article {KlappeneckerRottelerII02,
    AUTHOR = {Klappenecker, Andreas and R\"otteler, Martin},
     TITLE = {Beyond stabilizer codes. {II}. {C}lifford codes},
   JOURNAL = {IEEE Trans. Inform. Theory},
  FJOURNAL = {Institute of Electrical and Electronics Engineers.
              Transactions on Information Theory},
    VOLUME = {48},
      YEAR = {2002},
    NUMBER = {8},
     PAGES = {2396--2399},
}

@misc{KnillII96,
  title={Non-binary Unitary Error Bases and Quantum Codes}, 
  author={E. Knill},
  year={1996},
  eprint={quant-ph/9608048},
  archivePrefix={arXiv},
  primaryClass={quant-ph}, 
}

@misc{KnillI96,
  title={Group Representations, Error Bases and Quantum Codes}, 
  author={E. Knill},
  year={1996},
  eprint={quant-ph/9608049},
  archivePrefix={arXiv},
  primaryClass={quant-ph},
}

@article {KnillLaflamme97,
    AUTHOR = {Knill, Emanuel and Laflamme, Raymond},
     TITLE = {Theory of quantum error-correcting codes},
   JOURNAL = {Phys. Rev. A (3)},
  FJOURNAL = {Physical Review. A. Third Series},
    VOLUME = {55},
      YEAR = {1997},
    NUMBER = {2},
     PAGES = {900--911},
}

@book{Brown82,
series = {Graduate texts in mathematics},
volume = {87},
publisher = {Springer},
year = {1982},
title = {Cohomology of groups},
address = {New York},
author = {Brown, Kenneth S},
}

@incollection{CeccheriniSilbersteinTullio22,
series = {Springer Monographs in Mathematics},
publisher = {Springer International Publishing AG},
year = {2022},
title = {Representations of Finite Group Extensions via Projective Representations},
copyright = {The Author(s), under exclusive license to Springer Nature Switzerland AG 2022},
address = {Switzerland},
author = {Ceccherini-Silberstein, Tullio and Tolli, Filippo and Scarabotti, Fabio},
}

@article{Cheng15,
title = {A character theory for projective representations of finite groups},
journal = {Linear Algebra and its Applications},
volume = {469},
pages = {230-242},
year = {2015},
author = {Chuangxun Cheng}
}

@book{Waldron18,
series = {Applied and Numerical Harmonic Analysis},
publisher = {Springer Nature},
year = {2018},
title = {An Introduction to Finite Tight Frames},
edition = {1st Edition 2018},
copyright = {Springer Science+Business Media, LLC 2018},
address = {New York, NY},
author = {Waldron, Shayne F. D},
}

@phdthesis{Gottesman97,
  author = {Gottesman, Daniel Eric},
  school = {California Institute of Technology},
  title = {Stabilizer Codes and Quantum Error Correction},
  year = {1997},
}

@article{KribsLaflammePoulin2005,
  title = {Unified and Generalized Approach to Quantum Error Correction},
  author = {Kribs, David and Laflamme, Raymond and Poulin, David},
  journal = {Phys. Rev. Lett.},
  volume = {94},
  year = {2005},
}

@book{NielsenChuang2010,
  author = {Michael A. Nielsen and Isaac L. Chuang},
  edition = {10th Anniversary Edition},
  publisher = {Cambridge University Press},
  title = {Quantum Computation and Quantum Information},
  year = {2010}
}

@article{Connes1989,
	author = {Connes, Alain},
	journal = {Ergodic Theory Dynam. Systems},
	number = {2},
	pages = {207--220},
	title = {Compact metric spaces, {F}redholm modules, and hyperfiniteness},
	volume = {9},
	year = {1989}
}

@misc{Rieffel2025,
  title={Quantum Hamming Metrics}, 
  author={Marc A. Rieffel},
  year={2025},
  eprint={2507.23046},
  archivePrefix={arXiv},
  primaryClass={math.OA},
}

@book{ConnesNCGBook,
	author = {Connes, Alain},
	pages = {xiv+661},
	publisher = {Academic Press, Inc., San Diego, CA},
	title = {Noncommutative geometry},
	year = {1994}
}

@article{Kaad24,
	author = {Kaad, Jens},
	journal = {Pure Appl. Funct. Anal.},
	number = {5},
	pages = {1175--1211},
	title = {External products of spectral metric spaces},
	volume = {9},
	year = {2024}
}

@article {Latremoliere2015Dual,
    AUTHOR = {Latr\'emoli\`ere, Fr\'ed\'eric},
     TITLE = {The dual {G}romov-{H}ausdorff propinquity},
   JOURNAL = {J. Math. Pures Appl. (9)},
  FJOURNAL = {Journal de Math\'ematiques Pures et Appliqu\'ees. Neuvi\`eme
              S\'erie},
    VOLUME = {103},
      YEAR = {2015},
    NUMBER = {2},
     PAGES = {303--351},
      ISSN = {0021-7824,1776-3371},
   MRCLASS = {46L89 (46L30 58B34)},
  MRNUMBER = {3298361},
MRREVIEWER = {Efton\ Park},
       DOI = {10.1016/j.matpur.2014.04.006},
}

@article {Latremoliere2017TriangleIneq,
    AUTHOR = {Latr\'emoli\`ere, Fr\'ed\'eric},
     TITLE = {The triangle inequality and the dual {G}romov-{H}ausdorff
              propinquity},
   JOURNAL = {Indiana Univ. Math. J.},
  FJOURNAL = {Indiana University Mathematics Journal},
    VOLUME = {66},
      YEAR = {2017},
    NUMBER = {1},
     PAGES = {297--313},
      ISSN = {0022-2518,1943-5258},
   MRCLASS = {46L89 (46L35 58B34)},
  MRNUMBER = {3623412},
MRREVIEWER = {Hyun\ Ho\ Lee},
       DOI = {10.1512/iumj.2017.66.5954},
}

@article {AguilarLatremoliereRainone,
    AUTHOR = {Aguilar, Konrad and Latr\'emoli\`ere, Fr\'ed\'eric and
              Rainone, Timothy},
     TITLE = {Bunce-{D}eddens algebras as quantum {G}romov-{H}ausdorff
              distance limits of circle algebras},
   JOURNAL = {Integral Equations Operator Theory},
  FJOURNAL = {Integral Equations and Operator Theory},
    VOLUME = {94},
      YEAR = {2022},
    NUMBER = {1},
     PAGES = {Paper No. 2, 42},
      ISSN = {0378-620X,1420-8989},
   MRCLASS = {46L89 (46L30 58B34)},
  MRNUMBER = {4353489},
MRREVIEWER = {Soumalya\ Joardar},
       DOI = {10.1007/s00020-021-02678-w},
}

@article {LatremoliereCategories2019,
    AUTHOR = {Latr\'emoli\`ere, Fr\'ed\'eric},
     TITLE = {Actions of categories by {L}ipschitz morphisms on limits for
              the {G}romov-{H}ausdorff propinquity},
   JOURNAL = {J. Geom. Phys.},
  FJOURNAL = {Journal of Geometry and Physics},
    VOLUME = {146},
      YEAR = {2019},
     PAGES = {103481, 31},
      ISSN = {0393-0440,1879-1662},
   MRCLASS = {46L89 (18B99 46L30 58B34)},
  MRNUMBER = {3994972},
       DOI = {10.1016/j.geomphys.2019.103481},
}

@article {Gelander2013,
    AUTHOR = {Gelander, Tsachik},
     TITLE = {Limits of finite homogeneous metric spaces},
   JOURNAL = {Enseign. Math. (2)},
  FJOURNAL = {L'Enseignement Math\'ematique. Revue Internationale. 2e
              S\'erie},
    VOLUME = {59},
      YEAR = {2013},
    NUMBER = {1-2},
     PAGES = {195--206},
      ISSN = {0013-8584},
   MRCLASS = {54E35 (54E50)},
  MRNUMBER = {3113605},
MRREVIEWER = {David\ Gauld},
       DOI = {10.4171/LEM/59-1-8},
}

@article {ArzhantsevaCherix2024,
    AUTHOR = {Arzhantseva, Goulnara and Cherix, Pierre-Alain},
     TITLE = {Quantifying metric approximations of discrete groups},
   JOURNAL = {Ann. Fac. Sci. Toulouse Math. (6)},
  FJOURNAL = {Annales de la Facult\'e{} des Sciences de Toulouse.
              Math\'ematiques. S\'erie 6},
    VOLUME = {33},
      YEAR = {2024},
    NUMBER = {3},
     PAGES = {831--877},
      ISSN = {0240-2963,2258-7519},
   MRCLASS = {20E26 (03C20 20C99 20F69)},
  MRNUMBER = {4822917},
MRREVIEWER = {Egle\ Bettio},
       DOI = {10.5802/afst.1788},
}

@article {OlesenPedersenTakesaki1980,
    AUTHOR = {Olesen, Dorte and Pedersen, Gert K. and Takesaki, Masamichi},
     TITLE = {Ergodic actions of compact abelian groups},
   JOURNAL = {J. Operator Theory},
  FJOURNAL = {Journal of Operator Theory},
    VOLUME = {3},
      YEAR = {1980},
    NUMBER = {2},
     PAGES = {237--269},
      ISSN = {0379-4024},
   MRCLASS = {46L55 (22D40)},
  MRNUMBER = {578942},
MRREVIEWER = {H. Halpern},
}

@article {ArvesonPrice2003,
    AUTHOR = {Arveson, William and Price, Geoffrey},
     TITLE = {The structure of spin systems},
   JOURNAL = {Internat. J. Math.},
  FJOURNAL = {International Journal of Mathematics},
    VOLUME = {14},
      YEAR = {2003},
    NUMBER = {2},
     PAGES = {119--137},
      ISSN = {0129-167X,1793-6519},
   MRCLASS = {46L35 (46N50 47N50 82B10)},
  MRNUMBER = {1966768},
MRREVIEWER = {Marie\ Choda},
       DOI = {10.1142/S0129167X03001673},
}

@article {RainsNBQC99,
    AUTHOR = {Rains, Eric M.},
     TITLE = {Nonbinary quantum codes},
   JOURNAL = {IEEE Trans. Inform. Theory},
  FJOURNAL = {Institute of Electrical and Electronics Engineers.
              Transactions on Information Theory},
    VOLUME = {45},
      YEAR = {1999},
    NUMBER = {6},
     PAGES = {1827--1832},
      ISSN = {0018-9448,1557-9654},
   MRCLASS = {94B05 (94B65)},
  MRNUMBER = {1720636},
       DOI = {10.1109/18.782103},
}

@INPROCEEDINGS{WangZengGrasslSanders,
  author={Wang, Yun-Jiang and Zeng, Bei and Grassl, Markus and Sanders, Barry C.},
  booktitle={2013 IEEE International Symposium on Information Theory}, 
  title={Stabilizer formalism for generalized concatenated quantum codes}, 
  year={2013},
  volume={},
  number={},
  pages={529-533},
  doi={10.1109/ISIT.2013.6620282}}

@article {RieffelMetricsCompact98,
    AUTHOR = {Rieffel, Marc A.},
     TITLE = {Metrics on states from actions of compact groups},
   JOURNAL = {Doc. Math.},
  FJOURNAL = {Documenta Mathematica},
    VOLUME = {3},
      YEAR = {1998},
     PAGES = {215--229},
      ISSN = {1431-0635,1431-0643},
   MRCLASS = {46L87 (46L55 58B30 60B10)},
  MRNUMBER = {1647515},
MRREVIEWER = {A.\ I.\ Danilenko},
}

@article {RieffelMetrics99,
    AUTHOR = {Rieffel, Marc A.},
     TITLE = {Metrics on state spaces},
   JOURNAL = {Doc. Math.},
  FJOURNAL = {Documenta Mathematica},
    VOLUME = {4},
      YEAR = {1999},
     PAGES = {559--600},
      ISSN = {1431-0635,1431-0643},
   MRCLASS = {46L87 (58B34)},
  MRNUMBER = {1727499},
MRREVIEWER = {Jacek\ Brodzki},
}

@incollection {RieffelqGH04,
    AUTHOR = {Rieffel, Marc A.},
     TITLE = {Gromov-{H}ausdorff distance for quantum metric spaces},
      NOTE = {Appendix 1 by Hanfeng Li,
              Gromov-Hausdorff distance for quantum metric spaces. Matrix
              algebras converge to the sphere for quantum Gromov-Hausdorff
              distance},
   JOURNAL = {Mem. Amer. Math. Soc.},
  FJOURNAL = {Memoirs of the American Mathematical Society},
    VOLUME = {168},
      YEAR = {2004},
    NUMBER = {796},
     PAGES = {1--65},
      ISSN = {0065-9266,1947-6221},
   MRCLASS = {46L87 (53C23 58B34 60B10)},
  MRNUMBER = {2055927},
       DOI = {10.1090/memo/0796},
}

@article {ChristensenIvanAF06,
    AUTHOR = {Christensen, Erik and Ivan, Cristina},
     TITLE = {Spectral triples for {AF} {$C^*$}-algebras and metrics on the
              {C}antor set},
   JOURNAL = {J. Operator Theory},
  FJOURNAL = {Journal of Operator Theory},
    VOLUME = {56},
      YEAR = {2006},
    NUMBER = {1},
     PAGES = {17--46},
      ISSN = {0379-4024,1841-7744},
   MRCLASS = {46L87 (28A78 28A80 58B34 58J42)},
  MRNUMBER = {2261610},
MRREVIEWER = {J.\ D. Maitland Wright},
}

@misc{AAAKMqGH23,
Author = {Clay Adams and Konrad Aguilar and Esteban Ayala and Evelyne Knight and Chloe Marple},
Title = {Quantum Gromov-Hausdorff propinquity convergence of Christensen-Ivan quantum metrics on AF algebras},
Year = {2023},
Eprint = {arXiv:2312.16458},
}

@article {AguilarIndLim21,
    AUTHOR = {Aguilar, Konrad},
     TITLE = {Inductive limits of {$\rm C^*$}-algebras and compact quantum
              metric spaces},
   JOURNAL = {J. Aust. Math. Soc.},
  FJOURNAL = {Journal of the Australian Mathematical Society},
    VOLUME = {111},
      YEAR = {2021},
    NUMBER = {3},
     PAGES = {289--312},
      ISSN = {1446-7887,1446-8107},
   MRCLASS = {46L89 (46L30 58B34)},
  MRNUMBER = {4337940},
MRREVIEWER = {Paul\ Jolissaint},
       DOI = {10.1017/S1446788720000130},
}

@article {AguilarLatremoliere15,
    AUTHOR = {Aguilar, Konrad and Latr\'emoli\`ere, Fr\'ed\'eric},
     TITLE = {Quantum ultrametrics on {AF} algebras and the
              {G}romov-{H}ausdorff propinquity},
   JOURNAL = {Studia Math.},
  FJOURNAL = {Studia Mathematica},
    VOLUME = {231},
      YEAR = {2015},
    NUMBER = {2},
     PAGES = {149--193},
      ISSN = {0039-3223,1730-6337},
   MRCLASS = {46L89 (46L30 58B34)},
  MRNUMBER = {3465284},
MRREVIEWER = {\'Eric\ Leichtnam},
       DOI = {10.4064/sm227-2-5},
}

@article {AguilarFelltopAF19,
    AUTHOR = {Aguilar, Konrad},
     TITLE = {Fell topologies for {AF}-algebras and the quantum propinquity},
   JOURNAL = {J. Operator Theory},
  FJOURNAL = {Journal of Operator Theory},
    VOLUME = {82},
      YEAR = {2019},
    NUMBER = {2},
     PAGES = {469--514},
      ISSN = {0379-4024,1841-7744},
   MRCLASS = {46L89 (46L30 58B34)},
  MRNUMBER = {4015960},
MRREVIEWER = {Christopher\ Jack\ Bourne},
       DOI = {10.7900/jot},
}

@article {AGRKLstrongLeibniz24,
    AUTHOR = {Aguilar, Konrad and Garcia, Stephan Ramon and Kim, Elena and
              Latr\'emoli\`ere, Fr\'ed\'eric},
     TITLE = {The strongly {L}eibniz property and the {G}romov-{H}ausdorff
              propinquity},
   JOURNAL = {J. Math. Anal. Appl.},
  FJOURNAL = {Journal of Mathematical Analysis and Applications},
    VOLUME = {529},
      YEAR = {2024},
    NUMBER = {1},
     PAGES = {Paper No. 127572, 22},
      ISSN = {0022-247X,1096-0813},
   MRCLASS = {46L89 (46L30 58B34)},
  MRNUMBER = {4620742},
MRREVIEWER = {M\u ad\u alina\ Roxana\ Buneci},
       DOI = {10.1016/j.jmaa.2023.127572},
}

@article {AguilarYu24,
    AUTHOR = {Aguilar, Konrad and Yu, Jiahui},
     TITLE = {The {F}ell topology and the modular {G}romov-{H}ausdorff
              propinquity},
   JOURNAL = {Proc. Amer. Math. Soc.},
  FJOURNAL = {Proceedings of the American Mathematical Society},
    VOLUME = {152},
      YEAR = {2024},
    NUMBER = {4},
     PAGES = {1711--1724},
      ISSN = {0002-9939,1088-6826},
   MRCLASS = {46L89 (46L30 58B34)},
  MRNUMBER = {4709237},
MRREVIEWER = {Huichi\ Huang},
       DOI = {10.1090/proc/16669},
}

@article{BenyKempfKribs2007I,
  title = {Generalization of Quantum Error Correction via the Heisenberg Picture},
  author = {B{\'e}ny, C{\'e}dric and Kempf, Achim and Kribs, David W.},
  journal = {Phys. Rev. Lett.},
  volume = {98},
  pages = {paper number 100502},
  year = {2007},
}

@article{BenyKempfKribs2007II,
  title = {Quantum error correction of observables},
  author = {B{\'e}ny, C{\'e}dric and Kempf, Achim and Kribs, David W.},
  journal = {Phys. Rev. A},
  volume = {76},
  pages = {paper number 042303},
  year = {2007},
}

@article{KribsLaflammePoulin2006,
    author = {Kribs, David and Laflamme, Raymond and Poulin, David},
    title = {Operator quantum error correction},
    journal = {Quantum Inf. Comput.},
    volume = {6},
    pages = {383-399},
    year = {2006}
}

\end{document}